\documentclass{amsart}
\usepackage{amssymb,amsmath,amscd,amsthm}
\usepackage{amsfonts,epsfig,latexsym,graphicx,amssymb}
\usepackage{color}
\usepackage{enumerate}
\usepackage{todonotes}

\usepackage{hyperref}

\numberwithin{equation}{section}

\date{\today}

\newcommand{\Z}{{\mathbb Z}}
\newcommand{\R}{{\mathbb R}}
\newcommand{\Q}{{\mathbb Q}}
\newcommand{\C}{{\mathbb C}}
\newcommand{\N}{{\mathbb N}}

\newcommand{\CP}{\mathbb{CP}}

\newcommand{\ve}{\mathbf{e}}

\newcommand{\mF}{\mathcal{F}}

\newcommand{\per}{\mathrm{per}}
\newcommand{\anti}{\mathrm{anti}}
\newcommand{\Hf}{H}
\newcommand{\Hfp}{\Hf^{\per}}
\newcommand{\Hfa}{\Hf^{\anti}}
\newcommand{\pE}{E^{\per}}
\newcommand{\aE}{E^{\anti}}
\newcommand{\Vp}{\mV^{\per}}
\newcommand{\Va}{\mV^{\anti}}

\newcommand{\mV}{\mathcal{V}}
\newcommand{\mT}{\mathcal{T}}
\newcommand{\mI}{\mathcal{I}}

\newcommand{\Th}{\Theta}

\newcommand{\fa}{ {\frac{\alpha}{2}}}
\newcommand{\ta}{ {\theta_+}}

\newcommand{\tmu}{\widetilde\mu}
\newcommand{\tc}{\widetilde{c}}

\newcommand{\SL}{\mathrm{SL}}

\newcommand{\dist}{\mathop{\mathrm{dist}}\nolimits}
\newcommand{\Leb}{\mathop{\mathrm{Leb}}\nolimits}
\newcommand{\Ind}{1\hspace{-2.6pt}\mathrm{I}}

\newcommand{\Si}{\Sigma^{-}}
\newcommand{\Spl}{\Sigma^{+}}
\newcommand{\mR}{\mathcal{R}}

\newcommand{\mmu}{\mu^-}

\newcommand{\indsymb}{\ast}
\newcommand{\ba}{\alpha_{\indsymb}}
\newcommand{\bl}{\lambda_{\indsymb}}

\DeclareMathOperator{\tr}{tr}
\DeclareMathOperator{\SymDiff}{\Delta}

\newcommand{\Id}{\mathop{\mathrm{Id}}}

\newcommand{\Ta}{a}

 \graphicspath{ {pics/} } 

\newtheorem{theorem}{Theorem}[section]
\newtheorem*{theorem*}{Theorem}
\newtheorem{lemma}[theorem]{Lemma}
\newtheorem{prop}[theorem]{Proposition}
\newtheorem{coro}[theorem]{Corollary}

\newtheorem{example}[theorem]{Example}
\theoremstyle{definition}
\newtheorem{remark}[theorem]{Remark}
\newtheorem{question}[theorem]{Question}
\newtheorem*{defi}{Definition}

\newcommand{\eps}{{\varepsilon}}
\newcommand{\Sc}{\mathbb{S}^1}  %TODO: interchange!

\title[Continuity of the intersection spectrum of the AMO]{Generalized Aubry-Andr\'e formula \\ and continuity of the intersection spectrum \\ of the Almost Mathieu operator}

\author[A.\ Gorodetski]{Anton Gorodetski}

\address{Department of Mathematics, University of California, Irvine, CA~92697, USA}

\email{asgor@uci.edu}

\thanks{A.\ G.\ was supported in part by NSF grant DMS--2247966.}

\author[V. Kleptsyn]{Victor Kleptsyn}

\address{CNRS, University of Rennes, IRMAR, UMR 6625 du CNRS}

\email{victor.kleptsyn@univ-rennes.fr}

\thanks{V.K. was supported in part by Centre Henri Lebesgue (ANR-11-LABX-0020-01)}

\begin{document}

\begin{abstract}
We consider the spectrum of the Almost Mathieu operator (AMO) and show that the moments of the restriction of the Lebesgue measure to the intersection spectrum $\Leb|_{\Si_{\alpha,\lambda}}$ are polynomials in coupling $\lambda$ with coefficients that are trigonometric polynomials in frequency $\alpha$. The statement can be considered as a generalization of the Aubry-Andr\'e formula for the measure of the spectrum of AMO. As a corollary, we obtain that the restriction of the Lebesgue measure to the intersection spectrum that we denote by $\mmu_{\alpha, \lambda}$ depends continuously on the parameters (frequency $\alpha$  and coupling $\lambda$) in weak-* topology.  

Moreover, we prove that the dependence is not just continuous but analytic in $\lambda$ and $C^{\infty}$ in $\alpha$ in a sense that an integral of an analytic test function~$\varphi(x)$ with respect to $\mmu_{\alpha, \lambda}$ has the same kind of dependence. In particular, this implies that %As a corollary, we show that 
%As the analyticity of the test function  here is required only in a neighborhood of the union spectrum~$\Sigma^+_{\alpha,\lambda}$, this implies that
the Lebesgue measure of the part of the spectrum~$\Si_{\alpha,\lambda}$ that lies between two gaps depends analytically on the coupling constant $\lambda$ and $C^{\infty}$ on the frequency $\alpha$ in an open domain (away from the critical coupling $\lambda=1$) where these gaps do not bifurcate. 
\end{abstract}

\maketitle

\tableofcontents

%\newpage

%\subsection*{List of TODOs:}

%\mbox{}

%\listoftodos[]{}

%\mbox{}

%\newpage

\section{Introduction}
%\subsection{Preliminaries and questions}

In this paper we provide several new results on the Almost Mathieu operator, arguably one of the most heavily studied models in spectral theory of the last few decades. 

\subsection{History and background}

\begin{figure}\label{f.HB}
    \centering
\includegraphics[width=0.7\linewidth]{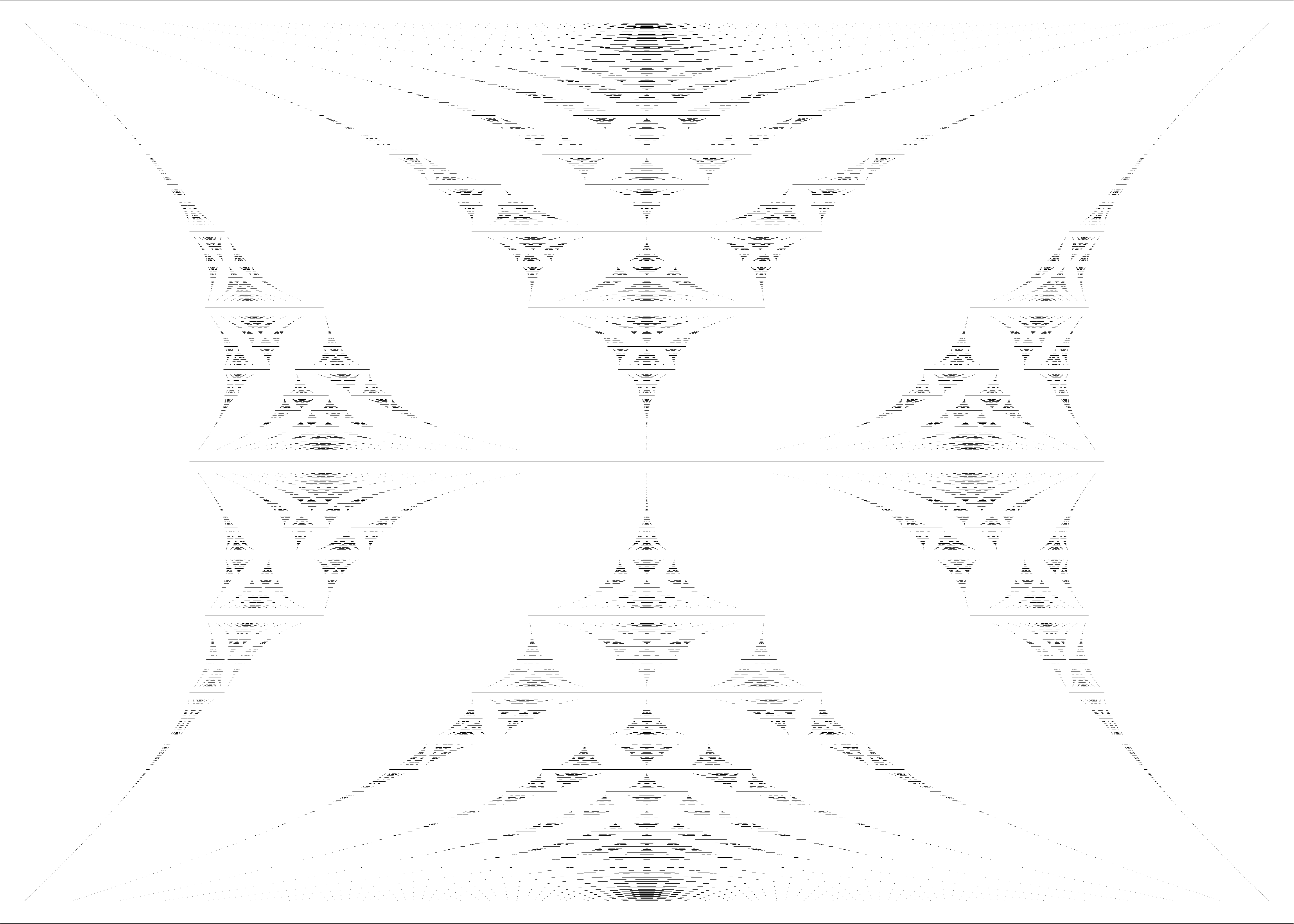}
    \caption{Hofstadter butterfly, $\lambda=1$}
    \label{fig:HB-critical}
\end{figure}

The \emph{Almost Mathieu operator} $H_{\alpha,\lambda,\theta}:\ell^2(\Z)\to \ell^2(\Z)$, defined by a frequency $\alpha\in [0,1)$, coupling constant $\lambda\in \R$, and a phase $\theta\in [0,1)$,
\begin{equation}\label{e.AM}
(H_{\alpha,\lambda,\theta} \psi)_n =  \psi_{n+1} + \psi_{n-1} + 2 \lambda  \cos 2\pi (n\alpha + \theta) \cdot \psi_n,
\end{equation}
is one of the most prominent models in the spectral theory of ergodic Schr\"odinger operators. The name is related to similarity to the Mathieu operator that was considered by Emile Mathieu in 1873 \cite{Math}. If $\lambda=1$, which is usually called a {\it critical case}, the operator (\ref{e.AM}) is sometimes called {\it Harper equation}; it was used by R.\,Peierls and his student R. G.\,Harper in the 1950s~\cite{Har} as a mathematical description of electrons on a 2D-lattice, acted on by a perpendicular homogeneous magnetic field. Lots of attention was attracted to the Almost Mathieu operator when in 1976 D.\,Hofstadter discovered the fractal structure of its spectrum \cite{Hof}. The {\it Hofstadter butterfly}  that visualizes a family of spectra of (\ref{e.AM}) for $\lambda=1$ parameterized by the frequency $\alpha\in [0,1)$, provided a particularly intriguing picture, see Fig.~\ref{f.HB}.
That started a period of increased interest in the spectral properties of operators with almost periodic potentials in general and in the Almost Mathieu operator in particular. In \cite{S82}, B.\,Simon characterized it as "the almost periodic flu''. 
%in the following way:{\it "In many years, flu sweeps the world. The actual strain varies from year to year; some years it has been Hong Kong flu, some years swine flu. In 1981, it was the almost periodic flu!"} 

As always for an ergodic Schr\"odinger operator, two main questions are about the spectral type of the operator and about the structure of the spectrum as a set. The spectral type of (\ref{e.AM}) was studied in \cite{AA, AD, AJ10, BLT, BLT85, CD,  Del, GJLS, J99, L1}, and the description obtained there is known as {\it Metal-Insulator Transition} due to the transition from a.c. spectrum for subcritical ($\lambda<1$) to singular continuous at critical ($\lambda=1$) and, depending on the arithmetic properties of frequency and phase, to p.p. or s.c. spectrum in supercritical ($\lambda>1$) regime. 
As for the structure of the spectrum of (\ref{e.AM}) as a set,  the two main conjectures were the Aubry-Andr\'e conjecture and the Ten Martini Problem. The Aubry-Andr\'e conjecture claims that the Lebesge measure of the spectrum of (\ref{e.AM}) for any irrational frequency is equal to $|4-4|\lambda||$; it was proven in a series of papers \cite{AK, AvMS, JK, L2}. The Ten Martini Problem is to show that the spectrum of (\ref{e.AM}) is a Cantor set for any $\alpha\not\in \mathbb{Q}$ and $\lambda\ne 0$. It was  formulated in \cite{Az, S82, S00} and resolved by A.\,Avila and S.\,Jitomirskaya in \cite{AJ}, following a series of partial results in \cite{AK, BS, CEY, HS, L2, Puig, Sin}. 

For a summary of the known results on the Almost Mathieu operator see~\cite[Section 9.12]{DF24}, and, more generally, for discrete Schr\"odinger operators with quasi-periodic potentials \cite[Chapter 9]{DF24}. For a survey of most recent results see~\cite{J22}.

\subsection{Continuity of the spectrum}\label{s:continuity}

We will denote by $\Sigma_{\alpha,\lambda, \theta}$ the spectrum of the Almost Mathieu operator (\ref{e.AM}). It is well known that for $\alpha\notin \Q$ the spectrum does not depend on the phase $\theta$, and can be denoted in this case by $\Sigma_{\alpha,\lambda}$.

If we do not assume that $\alpha\notin \Q$, then we can define the {\it intersection spectrum} $\bigcap_{\theta\in [0,2\pi)} \Sigma_{\alpha,\lambda, \theta}$, which is non-empty for $\lambda\in [0,1]$ and is known to be an empty set for $\lambda>1$ and any rational $\alpha$, see \cite[Corollary 3.2]{AvMS}.

Due to Aubry-Andr\'e duality, see \cite[Section 9.10]{DF24} for a modern exposition, for any $\alpha\notin \mathbb{Q}$ one has 
\[
    \Sigma_{\alpha,\lambda}=\lambda \Sigma_{\alpha,1/\lambda}.
\]
It is convenient to mimic this formula to extend the intersection spectrum to supercritical regime, and denote
\begin{equation}\label{eq.Si-def-super}
\Sigma_{\alpha,\lambda}^{-}=\left\{
  \begin{array}{ll}
    \bigcap_{\theta\in [0,2\pi)} \Sigma_{\alpha,\lambda, \theta}, & \hbox{if $\lambda\in [0,1]$;} \\
    \lambda\bigcap_{\theta\in [0,2\pi)} \Sigma_{\alpha,1/\lambda, \theta}, & \hbox{if $\lambda>1$.}
  \end{array}
\right.
\end{equation}
The {\it union spectrum} (sometimes also called {\it joint spectrum}) will be denoted by
$$
\Sigma_{\alpha,\lambda}^{+}=\bigcup_{\theta\in [0,2\pi)} \Sigma_{\alpha,\lambda, \theta}.
$$
Notice that for any irrational $\alpha$ and any$\lambda\in \mathbb{R}$ one has
$$
\Sigma_{\alpha,\lambda}^{+}=\Si_{\alpha,\lambda}=\Sigma_{\alpha,\lambda}.
$$
Also, by $|\Sigma_{\alpha,\lambda}^{-}|$ and $|\Sigma_{\alpha,\lambda}^{+}|$ we will denote the Lebesgue measure of each of the sets. 

Finally, by $\mmu_{\alpha, \lambda}$ we will denote the Lebesgue measure restricted to the set $\Sigma_{\alpha,\lambda}^{-}$. Notice that from \cite{AD} it follows that for any irrational $\alpha$ one has $\text{\rm supp}\,\mmu_{\alpha, \lambda} =\Sigma_{\alpha,\lambda}.$

Continuity (or the lack of thereof) of the defined objects and quantities, $\Sigma_{\alpha,\lambda}^{-}, \Sigma_{\alpha,\lambda}^{+}, |\Sigma_{\alpha,\lambda}^{-}|, |\Sigma_{\alpha,\lambda}^{+}|,$ and $\mmu_{\alpha, \lambda}$, attracted a lot of attention. Let us list some of the known results:

\begin{enumerate}
    \item  The set $\Sigma_{\alpha,\lambda}^{+}$ does depend continuously on the parameters in Hausdorff metric, see \cite[Theorem 3.7]{AS}. Moreover, it was shown that $\Sigma_{\alpha,\lambda}^{+}$ depends $\frac{1}{2}$-H\"older continuously  on the frequency $\alpha$, see \cite[Proposition 7.1]{AvMS}. It was mentioned in \cite{AvMS} that $\Sigma_{\alpha,\lambda}^{+}$ {\it ``...is the more natural object {\rm [than $\Sigma_{\alpha,\lambda}^{-}$]} in that it has a set
theoretic continuity in $\alpha$''}. 
    \item At the same time, it was noted in \cite[Section 7]{AvMS} that neither $\Sigma_{\alpha,\lambda, \theta}$  with a given fixed $\theta$ nor $|\Sigma_{\alpha,\lambda}^{+}|$ can depend continuously on $\alpha$ at the rational frequencies $\alpha$;
    \item The measure of the intersection spectrum is given by $|\Sigma_{\alpha,\lambda}^{-}|=|4-4\lambda|$, see \cite{AvMS, JK}, and hence depends continuously on the parameters. 

\item\label{i:converge} In \cite{JM12} it was shown, among other results, that in subcritical regime for any irrational $\alpha$ there exists a sequence of rationals $\frac{p_n}{q_n}\to \alpha$ such that the indicator functions converge:
$$
\Ind_{\Sigma_{\frac{p_n}{q_n}, \lambda}^+}\to \Ind_{\Sigma_{\alpha, \lambda}}\ \ \text{\rm and}\ \ \Ind_{\Sigma_{\frac{p_n}{q_n}, \lambda}^-}\to \Ind_{\Sigma_{\alpha, \lambda}}\ \text{\rm Lebesgue a.e.}
$$
\end{enumerate}

\begin{remark}
Many of the results mentioned here hold also in the case of general analytic \cite{JK, JM12, Sha}, smooth \cite{Z19}, or even just H\"older continuous quasiperiodic potentials \cite{JM14} in the operator (\ref{e.AM}). But in this paper we concentrate exclusively on the Almost Mathieu operator. 
\end{remark}

\begin{figure}
    \centering
    \includegraphics[width=0.8\linewidth]{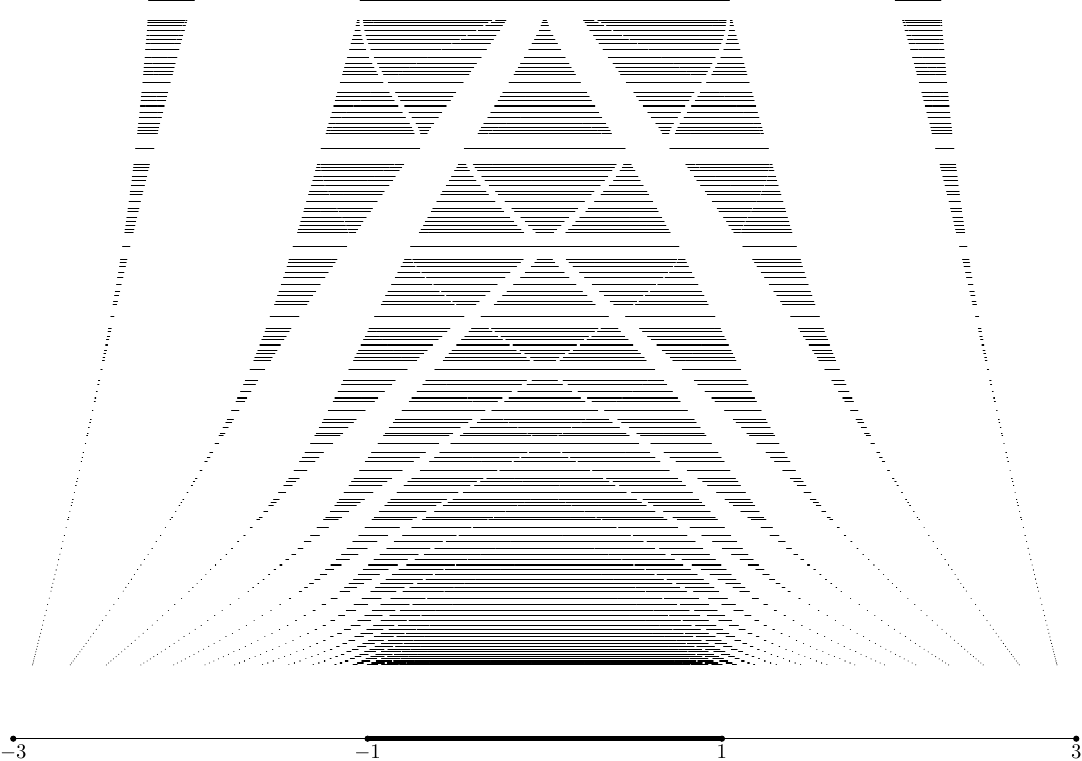}
    \caption{Hofstadter butterfly for $\lambda=1/2$, magnified near $\alpha=0$, with the intersection spectrum $\Si_{0,1/2}=[-1,1]$ drawn in bold.}
    \label{fig:HB-magnified}
\end{figure}

These facts paint somewhat inconsistent and even confusing picture. Our first main result deals with that inconsistency, and shows that not only the union spectra depend continuously on the parameters, but the intersection spectra as well, if one understand the continuity correctly, namely, as a continuity of the family of measures $\mmu_{\alpha, \lambda}$ (i.e. restrictions of the Lebesgue measure to $\Sigma_{\alpha, \lambda}^-$): 

\begin{theorem}\label{t.Leb-cont}
    The measures $\mmu_{\alpha, \lambda}$ depend continuously on $(\alpha,\lambda)$. %in the domain $\lambda<1$.
    In particular, for any $\frac{p}{q}\in\Q$ and any $\lambda\ne 1$, one has 
    \[
        \mmu_{\alpha, \lambda} \to \mmu_{\frac{p}{q}, \lambda}\quad \text{as} \quad \alpha\to \frac{p}{q}.
    \]
\end{theorem}
The following immediate corollary is illustrated by Fig.~\ref{fig:HB-magnified}: 
\begin{coro}\label{t.Leb-0}
For any $\lambda\in (0,1)$, as $\alpha\to 0$, the measures~$\mmu_{\alpha,\lambda}$ weakly converge to $\Leb|_{[-2(1-\lambda), 2(1-\lambda)]}$.
\end{coro}

%  It is known that for every irrational~$\alpha$, it does not depend on the phase $\theta_0$, so one can denote
% \[
% \Sigma_{\alpha,\lambda}=\sigma(H_{\alpha,\lambda;\theta_0}).
% \]

The proof of Theorem \ref{t.Leb-cont} is based on the properties of the moments of measures~$\mmu_{\alpha, \lambda}$, which leads us to our second result.

\subsection{Moments of the Lebesgue measure on the spectrum}

Our next result provides a  description of  moments of the measures $\mmu_{\alpha, \lambda}=\text{\rm Leb}|_{\Si_{\alpha, \lambda}}$ in subcritical regime. Namely, let us denote 
\begin{equation}\label{e.ck}
 c_{k}(\alpha,\lambda):=\int_{\mathbb{R}}x^k\, d\mu_{\alpha, \lambda}^-\equiv \int_{\Sigma_{\alpha,\lambda}^-} E^{k} \, dE.
\end{equation}
The odd moments $c_{2k+1}$ vanish due to the symmetry $E\mapsto -E$ of the spectrum of~\eqref{e.AM}, so we will consider only the even moments.

For $k=0$, due to Aubry-Andr\'e formula for the measure of the spectrum one has
$$
c_0(\alpha, \lambda)=|\Sigma_{\alpha,\lambda}^-|=|4-4\lambda|.
$$

Therefore, the following result can be considered as a generalized version of the Aubry-Andr\'e formula:

\begin{theorem}\label{t:poly}
For any $k\in \mathbb{N}$, and $\lambda\ge 0$, and any $\alpha\in \mathbb{R}$, the moment $c_{2k}$ can be represented as 
\begin{equation}\label{eq.c2k}    
c_{2k}=|P_{2k}(\alpha,\lambda)|=\left\{
  \begin{array}{ll}
    P_{2k}(\alpha,\lambda), & \hbox{if $\lambda\in [0,1];$}\\
    -P_{2k}(\alpha,\lambda), & \hbox{if $\lambda>1$,}
  \end{array}
\right.
\end{equation}
where $P_{2k}(\alpha,\lambda)$ is a polynomial of degree $2k+1$ as a function of $\lambda$ with the leading coefficient $-\frac{2^{2k+2}}{2k+1}$ and other coefficients that are trigonometric polynomials in $2\pi\alpha$, i.e. 
\begin{equation}\label{eq:P}
    P_{2k}(\alpha,\lambda) = \sum_{j,r\le 2k+1} q_{j,r} \cos^r 2\pi\alpha \cdot  \lambda^j.
\end{equation}
% such that for the moments $c_{2k}(\alpha,\lambda):=\int_{\mathbb{R}}x^{2k}\, d\mu_{\alpha, \lambda}^-\equiv \int_{\Sigma_{\alpha,\lambda}^-} E^{2k} \, dE$ we have
% \[
% \forall \wlambda\in[0,1] \quad \forall \alpha\in \R \quad c_{2k}(\alpha,\lambda) = P_{2k}(\alpha,\lambda).
% \]
% Moreover, this equality holds for all $\alpha$, including rational values, if one considers the intersection spectrum: 
% \begin{equation}\label{eq:int-c-2k}
%     \forall \lambda\in[0,1] \quad \forall \alpha\in \R \quad
%     \int_{\Si_{\alpha,\lambda}} E^{2k} \, dE = P_{2k}(\alpha,\lambda).
% \end{equation}
\end{theorem}

As we will see from the proof, polynomials $P_{2k}(\alpha,\lambda)$ can be calculated explicitly.
\begin{example}\label{ex:2-4}
    The second and the fourth moments in the subcritical regime $\lambda<1$ are given by 
    %\begin{multline*}        
    \[
        \frac{3}{16}c_2(\alpha,\lambda)= \frac{3}{16} \int_{\Sigma_{\alpha,\lambda}^-} E^2 \, dE 
         = 1-\frac{3}{2}\lambda (1+\cos 2\pi\alpha)+\frac{3}{2}\lambda^2 (1+\cos 2\pi\alpha) -  \lambda^3,
    \]
        % \\ = (2(1-\lambda))^3 - 12(\lambda-\lambda^2)(1-\cos \alpha),
    %\end{multline*}
    %\todo[inline]{Remove the second form for $c_2$, keep only the first one?}
    \begin{multline*}        
        \frac{5}{64}c_4(\alpha,\lambda)=\frac{5}{64} \int_{\Sigma_{\alpha,\lambda}^-} E^4 \, dE =
        \\  1 +\frac{5}{2}(\cos^3 2\pi\alpha +\cos^2 2\pi\alpha+\cos 2\pi\alpha+1) \lambda^2 + \frac{5}{4}(\cos 2\pi\alpha +1)^2 \lambda^4 
        \\ - \frac{5}{4}(\cos 2\pi\alpha +1)^2 \lambda
        - \frac{5}{2}(\cos^3 2\pi\alpha +\cos^2 2\pi\alpha+\cos 2\pi\alpha+1) \lambda^3
        - \lambda^5.
    \end{multline*}
\end{example}

We provide formulas for some other moments and discuss their properties in Section~\ref{s:conclusion}, see Example~\ref{ex:6-8}.

\begin{figure}
    \centering
    \includegraphics[width=0.9\linewidth]{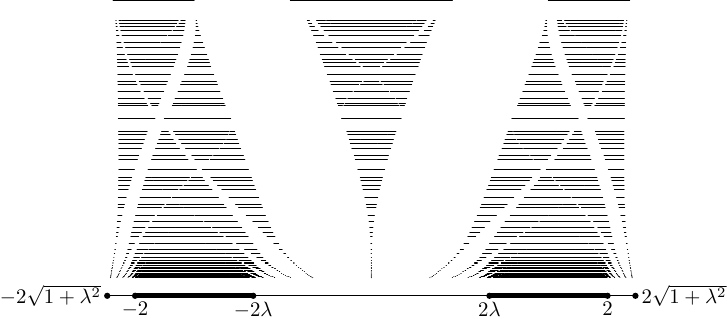}
    \caption{Hofstadter butterfly for a fixed $\lambda\in (0,1)$, magnified near $\alpha=\frac{1}{2}$, see Example \ref{ex:alphapi}.}
    \label{fig:HB-alphapi}
\end{figure}
The next example is illustrated by Fig. \ref{fig:HB-alphapi}.
\begin{example}\label{ex:alphapi}
For $\alpha=\frac{1}{2}$ and any $\lambda\in (0,1)$ one has 
$$
\Sigma_{\frac{1}{2}, \lambda}^+=[-2\sqrt{1+\lambda^2}, 2\sqrt{1+\lambda^2}], \ \ \Si_{\frac{1}{2},\lambda}=[-2, -2\lambda]\cup [2\lambda, 2],
$$
and
$$
 c_{2k}(\frac{1}{2}, \lambda)=-\frac{2^{2k+2}}{2k+1}(\lambda^{2k+1}-1), \ k\in \Z^+.
$$
\end{example}
As the reader may notice from Example~\ref{ex:2-4}, the polynomials $P_{2k}(\alpha,\lambda)$ also satisfy an Aubry-Andr\'e-type relation:
\begin{prop}\label{p.A-A-duality}
    For any $k\in \mathbb{Z}^+$ one has
\begin{equation}\label{eq:A-A-relation}
    P_{2k}(\alpha,\lambda)= -\lambda^{2k+1}P_{2k}\left(\alpha,\frac{1}{\lambda}\right).    
\end{equation}
\end{prop}

\subsection{Smoothness and analyticity of the Lebesgue measure of the spectrum}

Since the measures $\mu_{\alpha, \lambda}^-$ depend continuously on parameters, the $\mu_{\alpha, \lambda}^-$-measure of an interval with the end points in two gaps of the spectrum (i.e. the Lebesgue measure of the intersection spectrum between two given gaps)  must also depend continuously on parameters in a domain where the gaps persist. It is therefore natural to ask whether it has a better regularity then just being continuous. It turns out that it is actually analytic in $\lambda$ and $C^{\infty}$-smooth in $\alpha$, and this is our next result. %the last main result of this paper. 

\begin{figure}
    \centering
    \includegraphics[width=0.7\linewidth]{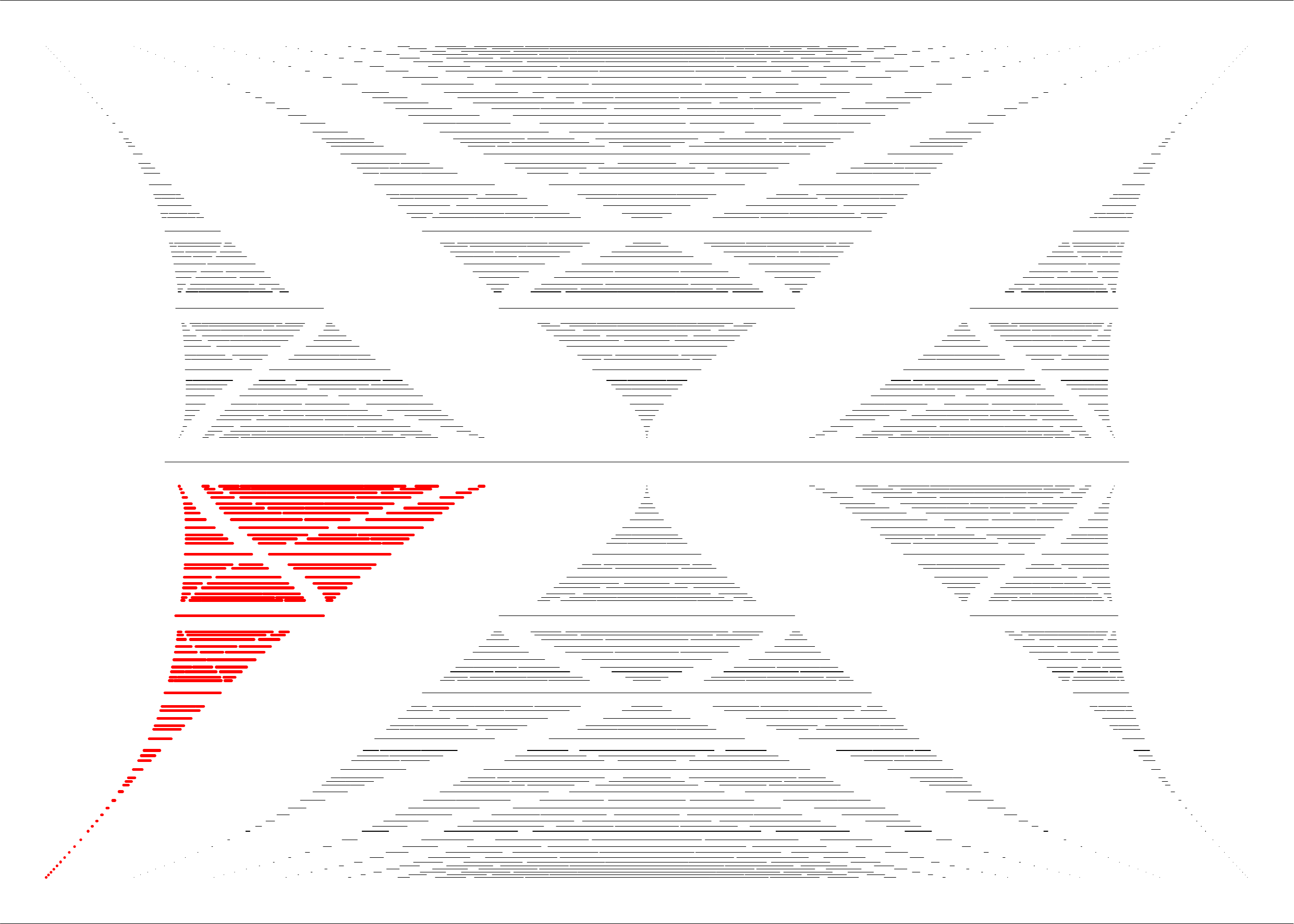}
    \caption{Hofstadter butterfly for $\lambda=1/2$ 
    and its part that for $0<\alpha<\frac{1}{2}$ lies between the gaps with labels $0$ and $1$ (highlighted).}
    \label{fig:HB-critical2}
\end{figure}

\begin{defi}
Take any $\alpha\notin \mathbb{Q}$ and any gap in $\Sigma_{\alpha, \lambda}$ (specified, for instance, by its label in the sense of the gap labeling theorem). Let $D$ be a domain in the $(\alpha,\lambda)$ plane, contained either entirely in the subcritical region $\lambda<1$, or in supercritical region $\lambda>1$.  
We say that the gap \emph{does not bifurcate} over a domain $D\subset \R\times \left((0,1)\cup (1, +\infty)\right)$ if its continuation (as a gap in $\Sigma_{\alpha, \lambda}^+$) is nondegenerate at every $(\alpha,\lambda)\in D$ (including the rational values of $\alpha$).
\end{defi}
 For instance, in Fig. \ref{fig:HB-critical2} one can see that the gap associated with the label  $m=1$ (one of the two largest gaps on the picture) bifurcates at $\alpha=\frac{1}{2}$.

Consider the part of $\Si_{\alpha,\lambda}$ that lies between two gaps (see Fig. \ref{fig:HB-critical2} for an example),  {and assume that these gaps} do not bifurcate over some domain~$D$. The following theorem holds.

\begin{theorem}\label{t.analytic}
The measure of the part of $\Sigma_{\alpha,\lambda}^-$ that lies between two gaps depends analytically on $\lambda$ and $C^{\infty}$-smoothly on $\alpha$ in any open domain (contained either entirely in the subcritical region $\lambda<1$, or in supercritical region $\lambda>1$) where these gaps do not bifurcate. 
\end{theorem}

\begin{remark}
If for some specific $\lambda\in (0,1)$ (or $\lambda>1$) there is a pair of gaps that do not bifurcate as the frequency changes within an interval $(\alpha_1, \alpha_2)$, then in Theorem \ref{t.analytic} one can take $D=(\alpha_1, \alpha_2)\times (0,1)$ (or $D=(\alpha_1, \alpha_2)\times (1, \infty)$, correspondingly). This is due  to the recent solution of the Dry Ten Martini Problem~\cite{AYZ} that claims that all gaps in the spectrum of the AMO with non-critical coupling that are prescribed by the Gap Labeling Theorem must actually be open. See \cite[Section 1.1]{AYZ} for a review of earlier results on the Dry Ten Martini Problem for the AMO. 
\end{remark}

\begin{remark}
It is interesting to take a look at the shape of a domain in $(\alpha, \lambda, E)$-space formed by continuations of a given gap in $\Sigma_{\alpha, \lambda}^+$.  For a fixed coupling, the boundary of the section of this domain in $(\alpha, E)$-plane is $\frac{1}{2}$-H\"older continuous due to \cite[Proposition 7.1]{AvMS}, see Fig.\ref{fig:HB-critical2}. At the same time, for a fixed frequency $\alpha$ it is reasonable to expect that the boundaries of the section in $(\lambda, E)$-plane are smooth or even analytic; indeed, in a related case of quasi-periodic  Hill's equation with Diophantine frequency analyticity of the gap boundaries was shown in \cite{PuigSimo}. More generally, it would be interesting to see whether a curve on $(\lambda, E)$-plane that corresponds to a given value of the integrated density of states (or, equivalently, to a given value of the rotation number of the corresponding Schr\"odinger cocycle) is a smooth or even analytic curve. One can compare these questions with similar questions on smoothness of boundaries of Arnold tongues \cite{Ar}, where the answers are affirmative \cite{GY, Sl2001}.     
\end{remark}

\subsection{Regularity of the families of measures} Theorem \ref{t.analytic} above claims that the Lebesgure measure of the part of the intersection spectrum depends analytically on $\lambda$ and $C^{\infty}$ on~$\alpha$. In fact, a much stronger analyticity statement regarding the family of measures $\mu_{\alpha, \lambda}^-$ holds. 
\begin{theorem}\label{t.anmu}
Suppose $\lambda_0\ne 1$,  $\alpha_0\in \mathbb{R}$, and an open set $U\subseteq \mathbb{R}$ are given such that $\Sigma_{\alpha_0, \lambda_0}^+\subset U$. Let $\varphi:U\to \R$ be an analytic function. Then there exists an open neighborhood  $D\subset \R^2$ of the point $(\alpha_0, \lambda_0)$ such that the integral
$$
\int \varphi(E) d\mu_{\alpha, \lambda}^-(E)
$$
depends analytically on $\lambda$ and $C^{\infty}$-smoothly on $\alpha$ in the domain~$D$.
\end{theorem}

\begin{remark}\label{r.particular}
In Theorem \ref{t.anmu}, one can also take an integral with respect to the measure $\mu_{\alpha, \lambda}^-$ restricted to a part of the spectrum between two gaps. Indeed, it suffices to take the set $U$ to be a union of two disjoint open sets $U=U_1\sqcup U_2$, where $U_1$ covers the part of the spectrum between the gaps, and  $U_2$ covers its complement, and extend the function $\varphi$ to be identically equal to zero in $U_2$. In particular, this observation turns Theorem \ref{t.analytic} into a partial case of Theorem \ref{t.anmu}.  
\end{remark}

%\subsection{Removing the singularity at $\lambda=1$}

Moreover, the singularity at the critical coupling $\lambda_0=1$ gets ``resolved'' when one passes from the family of measures $\mmu_{\alpha,\lambda}$ to the normalized family of \emph{probability} measures. Namely, the following statement holds:

\begin{theorem}\label{t.anmutilde}
There exists a family of probability measures $\{\tmu_{\alpha, \lambda}\}$, $\alpha\in \R$, $\lambda\ge 0$, such that the following holds:

\begin{enumerate}
    \item  For any $\lambda\ne 1$, we have $|4-4\lambda |\cdot \tmu_{\alpha, \lambda}=\mu_{\alpha, \lambda}^-$; 
    \item $\tmu_{\alpha,\lambda}$ depends on $(\alpha,\lambda)\in \Sc \times (0,\infty)$ continuously;
    \item For any $\alpha_0, \lambda_0$, an open set $U\subseteq \R$ with $\Sigma_{\alpha_0, \lambda_0}^+\subset U$, and an analytic function $\varphi:U\to \R$ there exists an open neighborhood  $D\subset \R^2$ of the point $(\alpha_0, \lambda_0)$ such that the integral
    \[
        \int \varphi(E) d\tmu_{\alpha, \lambda}(E)  
    \]
    is analytic in $\lambda$ and $C^{\infty}$-smooth in $\alpha$ in the domain~$D$.
\end{enumerate}

\end{theorem}

\begin{remark}

\begin{enumerate}[a)]
\item The family of measures $\mu_{\alpha, \lambda}^-$ degenerates to zero at the critical coupling. Theorem \ref{t.anmutilde} claims that normalizing these measures resolves this problem and makes the measures analytic in $\lambda$ everywhere, including the critical value of the coupling. As a simple example, let us apply Theorem \ref{t.anmutilde} to the function $\varphi\equiv 1$. Then we have
$$
\int \varphi(E) \, d\mmu_{\alpha, \lambda}(E)=|4-4\lambda|, \ \text{and} \ \ \int \varphi(x) \, d\tmu_{\alpha, \lambda}(x)=1. 
$$
Therefore, switching to normalized measures turned the non-differentiable at $\lambda=1$ function into an analytic one. 
\item In particular, this switching removes the sign change in the moments: for any~$k$, the moment 
\[
\tc_{2k}(\alpha,\lambda):=\int E^{2k}\, d\tmu_{\alpha,\lambda}(E)
\]
is equal to a function that is polynomial in $\lambda$ of degree $2k$ and is a trigonometric polynomial of~$\alpha$, and this equality holds for \emph{all}~$\alpha$ and $\lambda\in (0,+\infty)$. For instance, 
\[
\frac{3}{4}\tc_2(\alpha,\lambda)=1-\frac{1+  {3} \cos 2\pi \alpha}{2}\lambda + \lambda^2
\]
%\todo[inline]{Expression for $\tc_{4}$?}
(compare with Example~\ref{ex:2-4}).
\item For rational frequency $\alpha=\frac{p}{q}$ the measure  $\tmu_{\alpha, 1}$ is atomic, supported on $q$ points that form the set~$\Sigma_{\frac{p}{q},1}^{-}$.  {In particular, 
\[
\tmu_{0,1}=\delta_0,\quad \tmu_{\frac{1}{2},1}=\frac{1}{2}\left(\delta_{2}+\delta_{-2}\right)
\]
as one can see from Corollary~\ref{t.Leb-0} and Example~\ref{ex:alphapi} respectively. In general, the measure $\tmu_{\frac{p}{q},1}$ is equal to 
\[
q\cdot  \sum_{j=1}^q \frac{1}{|Q'_1(E_j)|} \delta_{E_j},
\]
where $Q_1(E)$ is given by Chambers' formula (see Eq.~\eqref{eq:Chambers} below), and $E_j$, $j=1,\dots, q$ are its roots. To verify this, it suffices to pass to the limit as $\lambda\to 1$ in~\eqref{eq:Si-Q}, see Remark~\ref{r:Si-Q} below. In particular, this implies that $\sum_{j=1}^q \frac{1}{|Q'_1(E_j)|}=q,$ which was initially observed in \cite{LW}. 
\item The family of measures $\tmu_{\alpha, \lambda}$ seems to be a very natural object. It would be interesting to understand whether it has any connection to the spectral properties of the initial AMO.
}
%\todo{What is the sum of these derivatives? Residues for the summation?}
\end{enumerate}
\end{remark}

\subsection{Main ideas and scheme of the proof}

 We start by generalizing in Section~\ref{s:rational} the approach of~\cite{AvMS} for the integrals of continuous functions. Namely, assume that $\alpha=\frac{p}{q}$, we are in the subcritical regime $\lambda<1$, and that we are integrating a function~$\varphi(E)=\Phi'(E)$. We use Chambers' formula to write the integral 
 \begin{equation}\label{eq:int-phi-scheme}
    \int \varphi(E) \, d\mmu_{\alpha,\lambda} 
 \end{equation}
 as an alternating sum~\eqref{eq:int-phi} of values $\Phi(E)$, taken over eigenvalues of the restrictions of the operator $H$ on some well-chosen finite-dimensional subspaces. It is done in Section~\ref{s:Chambers}--\ref{s:AvMS}, see Corollary~\ref{c.Phi}. Next, as the corresponding eigenfunctions are alternatively even and odd with respect to corresponding reflections, we transform this alternating sum into a difference of two traces of compositions of the form $\Phi(H)R$, where $H$ is a restriction of the Schr\"odinger operator to the corresponding finite-dimensional subspace, and $R$ is the  reflection, coming for the corresponding choice of the phase from the potential-preserving reflection $\mR: x\mapsto -x$ on the circle of phases (see Fig.~\ref{fig:circle-alpha}). This is done in Section~\ref{s:two}, see Proposition~\ref{p:Phi-tr}. We also note that as the reflection~$R$ interchanges the coordinates, each of these traces actually becomes the sum of elements of~$\Phi(H)$ along the corresponding anti-diagonal. Finally, in Section~\ref{s:A-A} we obtain a version of Aubry--Andr\'e duality in the case of the rational frequency, Proposition~\ref{p:Phi-tr-super}, showing that in the supercritical case $\lambda>1$ for the polynomials for the moments $c_{2k}$ one should simply change the sign of the formula.

 Section~\ref{s:infinite} is devoted to the passage to the language of infinite-dimensional operators. Namely, in Section~\ref{s:c-moments} we study further the traces that occur for the moments problem, and notice that each of these traces can be split into two ``traces'' of analogous infinite-dimensional matrices. These operators are not trace type, but have only finitely many nonzero elements on the main diagonal, thus we formally define their ``trace'' to be sum of these elements. One thus gets a sum of four ``traces'', see Proposition~\ref{p:moment-trace}, and the phases occurring in this formula correspond to the orbits of rotation by $\alpha$ that are fixed by the symmetry $x\mapsto -x$ acting on the circle. We then note that such a formula can be generalized to a wider class of functions~$\Phi$. Namely, in Section~\ref{s:decrease} we obtain an analogue of Combes--Thomas estimates: we show that for a function $\Phi$ that is analytic in a neighbourhood of the union spectrum $\Spl_{\alpha,\lambda}$, the elements of $\Phi(H)$ decay exponentially fast with the distance to the diagonal; see Proposition~\ref{p:Phi-bounds}. This allows us to define ``traces'' of their products $\Phi(H) R$ with the reflection matrices $R$ and to claim the ``four-traces'' formula in Section~\ref{s:four-traces}, see Theorem~\ref{t:four-traces}. We then complete proofs in Section~\ref{s:proofs}.

 Finally, in Section~\ref{s:ex-smooth} we address the question of the regularity in $\alpha$. Namely, in Theorem~\ref{t:four-traces} we claim that the result is $C^{\infty}$ in $\alpha$, but not analytic. We present an example that motivates our belief in absence of analyticity in $\alpha$ of the integral $\int \varphi(E) \, d\mmu_{\alpha,\lambda}(E)$ for a general analytic function~$\varphi(E)$.

 We discuss some open questions and possible generalisations in Section~\ref{s:conclusion}.
% \vspace{1cm}

% \todo[inline]{Rewrite the plan below}

% Next, considering the function $\Phi(E)=E^{2k+1}$ (corresponding to the moments $c_{2k}$), we note that 
%  in Section~\ref{s:infinite}

%  Next, assume that the function~$\varphi$ and hence its anti-derivative $\Phi$ are analytic in some neighbourhood of the spectrum $\Sigma^+_{\alpha,\lambda}$. Then, instead of a finite-dimensional space we consider the full action of $H_{\alpha,\lambda,\theta}$ on the space $\ell^2(\Z)$. On one hand, it turns out that the anti-diagonal elements of $\Phi(H_{\alpha,\lambda,\theta})$ exponentially decrease (what is a manifestation of Combes--Thomas estimates), and at the same time, they are analytic in $\lambda$ and $C^{\infty}$-smooth in~$\alpha$ (Proposition~\ref{???}). On the other hand,each of the finite-dimensional traces of $\Phi(H)R$ appearing in~\eqref{eq:int-phi} can be written as a sum or a difference of two associated ``traces'' for $\Phi(H_{\alpha,\lambda,\theta}) R_{\theta}$, where again the reflection $R_{\theta}$ comes from the restriction of the reflection $\mR$ on the orbit of the phase~$\theta$. It is done in Lemma~\ref{???}. As a corollary, we obtain formula~\eqref{eq:Phi-traces-int} for the integral~\eqref{eq:int-phi-scheme} that is a sum of four likewise ``traces'', each corresponding to one of four phases $\{0,\fa,\frac{1}{2},\frac{1+\alpha}{2}\}$, for which their orbits are preserved by the reflection~$\mR$ (see Fig.~\ref{fig:circle-alpha}). 
 
%  Next, we note that the right hand side of~\eqref{eq:Phi-traces-int} depends continuously on~$\alpha$.

 %\newpage

 \section{Finite dimensional operators: rational frequency, subcritical regime}\label{s:rational}

 Throughout this section, we will assume that $\alpha=\frac{p}{q}$ is rational and, until Section~\ref{s:A-A}, that we are in the subcritical case~$\lambda<1$. 
 
 \subsection{Intersection spectrum and Chambers' formula}\label{s:Chambers}

 We start by recalling the Chambers' formula. To do so, consider the Schr\"odinger cocycle, associated to the Almost Mathieu operator and the equation $H_{\alpha,\lambda,\theta}\psi=E\psi$ for its generalized eigenfunctions. Namely, let
\begin{equation}\label{eq:Pi-def}
    \Pi_{E,\theta} = \begin{pmatrix}
        E-V(\theta) & -1 \\
        1 & 0
    \end{pmatrix}, 
    \quad 
    T_{n,\theta}(E)=\Pi_{E,\theta+(n-1)\alpha}\dots \Pi_{E,\theta+\alpha} \Pi_{E,\theta}    
\end{equation}
(both $\Pi$ and $T$ depend on $\alpha$ and $\lambda$, but we do not indicate this dependence explicitly).
%\todo{Perhaps we should...}  %
Then $E$ belongs to the spectrum $\Sigma_{p/q,\lambda,\theta}$ if and only if the corresponding transfer matrix over the period $T_{q,\theta}$ is elliptic or parabolic; in other words, if and only if $|\tr T_{q,\theta}|\le 2$.

Now, consider the dependence of the trace $\tr T_{q,\theta}$ on the phase~$\theta$. This dependence is given by the Chambers formula (see~\cite{Cham} for the original paper by W.~Chambers, where it appeared initially, and \cite[Proposition~3.1]{AvMS} or \cite[Theorem~7.6.2]{DF22} for a short proof). This formula states that for $\alpha=\frac{p}{q}$ the trace of the transition operator along the period $T_{q,\theta}$ is equal to 
\begin{equation}\label{eq:Chambers}
\tr T_{q,\theta}(E) = - 2\lambda^q \cos 2\pi q\theta + Q_{\lambda}(E),    
\end{equation}
where $Q_{\lambda}(E)$ is a polynomial in $E$ and $\lambda$ of degree~$q$. In particular, for a given energy $E$ the trace takes its minimal value for the phase $\theta=0$ (as well as for all other $\theta=\frac{2\pi j}{q}$), and maximal value for all $\theta=\frac{(2j-1)\pi}{q}$, $j=1,\dots,q$. We choose an element of the latter family, fixing
\begin{equation}\label{eq:ta-def}
\ta=\begin{cases}
    \frac{\alpha}{2}, & p \text{ odd},\\
    \frac{1+\alpha}{2}, & p \text{ even};
\end{cases}    
\end{equation}
then, in both cases, one has 
\begin{equation}\label{eq:ta-choice}
    \ta-\alpha \equiv -\ta \mod 1, 
\end{equation}
a property that will be heavily used below.

 \subsection{Avron-von Mouche-Simon approach: even and odd eigenfunctions}\label{s:AvMS}

Following~\cite{AvMS}, note that for any phase $\theta$, the solutions $E$ to the equation $\tr T_{q,\theta}(E)=\pm 2$ correspond to the transfer matrices over the period with the eigenvalues $\pm 1$, and hence to the operator $H_{p/q,\lambda,\theta}$ possessing a $q$-periodic or $q$-antiperiodic generalized eigenfunction respectively: $\psi_{n+q}=\pm \psi_n$ for all~$n$. 

Consider the spaces $\Vp_q$ and $\Va_q$ of $q$-periodic and $q$-antiperiodic  functions~$\psi\in\ell^{\infty}(\Z)$ respectively. These spaces are $q$-dimensional, and for any phase $\theta$ are preserved by the corresponding Hamiltonian $H_{\frac{p}{q},\lambda,\theta}$ (as it commutes with a shift by $q$ due to $q\alpha=q\frac{p}{q}\in \Z$, and hence preserves its eigenspaces). Denote the restrictions of $H_{\frac{p}{q},\lambda,\theta}$ to these spaces by 
\[
\Hfp_{q,\theta}=\left.H_{\frac{p}{q},\lambda,\theta}\right|_{\mV_q^{\per}},
\quad
\Hfa_{q,\theta}=\left. H_{\frac{p}{q},\lambda,\theta}\right|_{\mV_q^{\anti}}.
\]
Due to the above arguments, their eigenvalues are the solutions to the corresponding equations $\tr T_{q,\theta}=\pm 2$. Denote these eigenvalues by $\pE_{j,\theta}$ and $\aE_{j,\theta}$ respectively: 
\[
\tr T_{q,\theta}(\pE_{j,\theta})=2, \quad \tr T_{q,\theta}(\aE_{j,\theta})=-2,
\]
ordering these eigenvalues in decreasing order:
\[
\pE_{1,\theta} \ge \dots\ge \pE_{q,\theta}, \quad 
\aE_{1,\theta} \ge \dots\ge \aE_{q,\theta}.
\]
Then, Chambers' formula immediately implies the following description for the intersection spectrum.
\begin{prop}\label{p:Si}
    The intersection spectrum for $\alpha=\frac{p}{q}$ 
    and $\lambda<1$ is given by
\begin{equation}\label{eq:Si-J}
    \Si_{\frac{p}{q},\lambda} = \bigcup_{j} J_j, \quad \text{where} \quad J_{j}=
    \begin{cases}
        [\aE_{j,0},\pE_{j,\ta}], & j \,\, \text{odd},\\
        [\pE_{j,\ta},\aE_{j,0}], & j \,\, \text{even}.    
    \end{cases}
\end{equation}
\end{prop}
\begin{remark}\label{r:Si-Q}
    Chambers' formula~\eqref{eq:Chambers} also allows to write~\eqref{eq:Si-J} in terms of the corresponding polynomial $Q_{\lambda}$ as a preimage
    \begin{equation}\label{eq:Si-Q}
        \Si_{\frac{p}{q},\lambda}= Q_{\lambda}^{-1}\left([-2(1-\lambda^q),2(1-\lambda^q)]\right)        
    \end{equation}
\end{remark}

\begin{proof}[Proof of Proposition~\ref{p:Si}]
    Due to Chambers' formula, for any phase $\theta$ the graph of $\tr T_{q,\theta}$ lies between the graphs of $\tr T_{q,0}$ and of $\tr T_{q,\ta}$ (see Fig.~\ref{fig:traces}). In particular, energy $E$ belongs to $\Si_{\alpha,\lambda}$ if and only if 
    \[
        \forall \theta \quad \tr T_{q,\theta}(E) \in [-2,2],
    \]
    that is equivalent to 
    \[
        -2\le \tr T_{q,0}(E), \text{ and } \tr T_{q,\ta}(E)\le 2;
    \]
    these inequalities define exactly the union of intervals~$J_j$.
\end{proof}

\begin{figure}
    \centering
    \includegraphics[width=0.7\textwidth]{Graph-10.pdf}
    \caption{Traces of the transfer matrices $T_{q,\theta}$ associated to the phases $\theta=0$ and $\theta=\ta$, drawn here for $q=3$. The periodic and antiperiodic eigenvalues $\pE_{j,\ta}, \aE_{j,0}$ (corresponding to the values of the trace~$2$ and~$-2$ respectively) are marked; the marking is a filled and empty circle for even and odd eigenvectors respectively. The intersection spectrum $\Si_{p/q,\lambda}=\bigcup_j J_j$ is marked by bold lines on the~$E$ axis below.} 
    \label{fig:traces}
\end{figure}

\begin{coro}\label{c.Phi}
    For rational $\alpha=\frac{p}{q}$, $\lambda<1$, and any continuous function $\varphi(E)=\Phi'(E)$ on the intersection spectrum $\Si_{\alpha,\lambda}$, one has 
    \begin{equation}\label{eq:int-phi}
        \int \varphi(E) \, d\mmu_{\alpha,\lambda}(E) = \sum_{j=1}^q (-1)^{j+1} \left(\Phi(\pE_{j,\ta}) - \Phi(\aE_{j,0})\right).        
    \end{equation}
\end{coro}
\begin{proof}
    Indeed, each of the summands in the right hand side of~\eqref{eq:int-phi} equals to the integral of~$\varphi$ over the corresponding~$J_j$, that is, to the increment of~$\Phi$ over~$J_j$. The~$(-1)^j$ coefficient comes from the order of the endpoints of~$J_j$ (see~\eqref{eq:Si-J}).
\end{proof}

 \subsection{Two traces formula}\label{s:two}

Now, the key idea of the work~\cite{AvMS} is to study which of the corresponding eigenfunctions are even or odd. Namely, the operator $\Hfa_{q,0}$ commutes with the reflection 
\begin{equation}\label{eq:R0}
R_0:(\psi_n)\mapsto (\psi_{-n}).     
\end{equation}
Hence, each of the eigenfunctions corresponding to (simple) eigenvalues~$\aE_{j,0}$ is either even or odd. In the same way, the operator $\Hfp_{q,\ta}$ commutes with the reflection 
\begin{equation}\label{eq:Rfa}
R_{\ta}:(\psi_n)\mapsto (\psi_{-1-n}),
\end{equation}
coming from the reflection $\theta\mapsto -\theta$ of the circle that preserves the orbit of the phase~$\ta$: from~\eqref{eq:ta-choice} one has
\[
  -(\ta+n\alpha) = \ta + (-1-n)\alpha.
\]
Hence, each of the eigenvalues $\Hfp_{q,\ta}$ is also either even or odd with respect to this reflection. Now, Theorems~4a and~4b from~\cite{AvMS} state that both sequences of eigenvalues appearing in the right hand side of~\eqref{eq:int-phi} alternate between even and odd. 
\begin{theorem}[J.~Avron, P.~H.~M.~van Mouche, B.~Simon, {\cite[Theorem~4a,~4b]{AvMS}}]\label{t.alt}
%\begin{itemize}
%    \item 
    The eigenvalues 
    \begin{equation}\label{eq:aE-0}
        \aE_{1,0}>\dots>\aE_{q,0}    
    \end{equation}
    of the operator $\Hfa_{q,0}$ alternately correspond to even and odd functions. In the same way,
%    \item 
    the eigenvalues 
    \begin{equation}\label{eq:pE-a-2}
    \pE_{1,\ta}>\dots>\pE_{q,\ta}
    \end{equation}
    of the operator $\Hfp_{q,\ta}$ alternately correspond to even and odd functions.
%\end{itemize}
\end{theorem}
For reader's convenience, we provide here a sketch of the proof (referring the reader to~\cite[Sections~4--5]{AvMS} for details).
\begin{proof}[Sketch of the proof]
Note that the listed eigenvalues stay simple: the inequality signs in~\eqref{eq:aE-0},~\eqref{eq:pE-a-2} are strict for any $\lambda>0$. Indeed, for instance $\aE_{1,0}$ and $\aE_{2,0}$ cannot merge as there is an intersection $\aE_{1,\ta}$ of the graph of $\tr T_{q,\ta}$ with the level $-2$ between them. This intersection cannot disappear as it corresponds to an eigenvalue of a self-adjoint operator $\Hfa_{q,\ta}$ in $q$-dimensional space. 

In the same way, $\aE_{2,0}$ and $\aE_{3,0}$ cannot merge because $\pE_{2,0}$ is between them. Again, this solution of $\tr T_{q,0}(E)=2$ cannot disappear as it corresponds to an eigenvalue of a self-adjoint operator $\Hfa_{q,0}$ in $q$-dimensional space. 
The same arguments apply to any other pair of consecutive eigenvalues, listed in~\eqref{eq:aE-0},~\eqref{eq:pE-a-2}, hence ensuring their simplicity.

As the eigenvalues listed in the conclusion of the theorem stay simple for any positive and finite $\lambda$, they cannot change their types (even/odd). Hence it suffices to establish the conclusion for just one value of the coupling constant~$\lambda$. This is done via a perturbation theory near $\lambda=\infty$, see~\cite[Sections~4--5]{AvMS} for details.
\end{proof}
\begin{remark}
    For other values of the phase, eigenvalues are not necessarily simple. For instance, $E=0$ is a double eigenvalue for~$\Hfa_{2,\frac{1}{4}}$ (that is, $\alpha=\frac{p}{q}=\frac{1}{2}$, and the phase $\theta=\frac{1}{4}$). We refer the reader to~\cite[Theorem~3]{BS} for a complete description of closed gaps. 
\end{remark}

\begin{remark}
    The conclusion of Theorem~\ref{t.alt}, stating that the eigenfunctions alternate between being even and odd, does not hold for other potentials~--- even for trigonometric polynomials of the form 
    \begin{equation}\label{eq:potential}
        V(\theta)=2\lambda_1\cos 2\pi \theta + 2\lambda_2\cos 4 \pi \theta.        
    \end{equation}
    For instance, take 
    \[
        \lambda_1=0.3, \quad \lambda_2=\frac{1}{5} \lambda_1=0.06, \quad \frac{p}{q}=\frac{2}{5}.
    \] 
    Then the eigenvalues $\pE_{j,\ta}$ of $\Hfp_{q,\ta}$ are 
    \[
      \{\underline{-1.727..}, -1.634.., 0.337.., \underline{0.967..},   2.057..\},
    \]
    where the underlined two eigenvalues correspond to the odd functions.
    %%%  Verification files: 
    %%%  Python: AMO.ipynb
    %%%  SciLab: perturbed-AM-pi.sce
\end{remark}

Theorem~\ref{t.alt} has allowed to Avron, van Mouche and Simon \cite{AvMS} to calculate the Lebesgue measure of the intersection spectrum. Namely, for 
\begin{multline}\label{eq:AvMS-Leb}
\left|\Si_{\frac{p}{q},\lambda}\right| = \int 1 \, d\mmu_{\frac{p}{q},\lambda} = \sum_{j=1}^q (-1)^{j+1} 
\left(\pE_{j,\ta} - \aE_{j,0}\right) =
\\
= \sum_{j\equiv 1 \mod 2}  \pE_{j,\ta} 
- \sum_{j\equiv 0 \mod 2} \pE_{j,\ta} 
+ \sum_{j\equiv 0 \mod 2} \aE_{j,0} 
- \sum_{j\equiv 1 \mod 2} \aE_{j,0},
\end{multline}
each of the four sums in the right hand side is the trace of the restriction of $\Hfa_{q,0}$ or $\Hfp_{q,\ta}$ on the subspace of even or odd functions. Then, by choosing a base in each of these spaces and explicitly calculating these traces, the authors of~\cite{AvMS} have performed the computation (see ~\cite[Proposition~3.3]{AvMS}), giving 
\[
\left|\Si_{\frac{p}{q},\lambda}\right| = 4(1-\lambda).
\]

Generalising their approach, we note that it is possible to count the alternating sums of eigenvalues in the right hand side of~\eqref{eq:int-phi} directly. Namely, the following proposition holds.

%c.Phi
\begin{prop}[Two traces formula]\label{p:Phi-tr}
    For any $\alpha=\frac{p}{q}$, any $\lambda<1$ and any continuous function $\Phi$ on $\Si_{\alpha,\lambda}$, the right hand side of~\eqref{eq:int-phi} can be written as
    \begin{equation}\label{eq:Phi-trace}
        \sum_{j=1}^q (-1)^{j+1} \left(\Phi(\pE_{j,\ta}) - \Phi(\aE_{j,0})\right) 
        = \tr \left(\Phi(\Hfp_{q,\ta}) R_{\ta}\right)
         -\tr \left(\Phi(\Hfa_{q,0}) R_0\right),
    \end{equation}
    where $R_0$, $R_{\ta}$ are the reflections, given by~\eqref{eq:R0} and~\eqref{eq:Rfa} respectively.

    In particular, if $\varphi=\Phi'$ exists and is continuous on $\Si_{\alpha,\lambda}$, then the integral of $\varphi$ over $\Si_{\alpha,\lambda}$ is equal to the right hand side of~\eqref{eq:Phi-trace}.
\end{prop}

\begin{proof}
    Indeed, the operators $\Hfa_{q,0}$ and $R_0$ commute, and their common eigenfunctions $\psi_{j,0}^{\anti}$ satisfy
    \[
        \Hfa_{q,0} \psi^{\anti}_{j,0} = \aE_{j,0}\psi^{\anti}_{j,0}, \quad R_{0} \psi^{\anti}_{j,0} = (-1)^{j+1} \psi^{\anti}_{j,0},
    \]
    where the sign is implied by the fact that the eigenfunctions are alternatively even and odd.
    Hence, 
    \begin{equation}\label{eq:Phi-anti}
        \sum_{j=1}^q (-1)^{j+1} \Phi(\aE_{j,0}) = \tr \left(\Phi(\Hfa_{q,0}) R_0\right),        
    \end{equation}
    as the $(-1)^{j+1}$ signs in the left hand side exactly correspond to the eigenvalues of~$R_0$.
    In the same way,
    \begin{equation}\label{eq:Phi-per}
        \sum_{j=1}^q (-1)^{j+1} \Phi(\pE_{j,\ta}) = \tr \left(\Phi(\Hfp_{q,\ta}) R_{\ta}\right),       
    \end{equation}
    and subtracting~\eqref{eq:Phi-anti} from~\eqref{eq:Phi-per}, we obtain the desired~\eqref{eq:Phi-trace}.
\end{proof}

%\newpage

\subsection{Aubry--André duality for rational frequencies}\label{s:A-A}

Let now $\alpha=\frac{p}{q}$, $\lambda>1$. It turns out that the formula for the integral from Proposition~\ref{p:Phi-tr} can be generalized to the supercritical regime with only a sign change: 

\begin{prop}\label{p:Phi-tr-super}
    For any rational $\alpha=\frac{p}{q}$, any $\lambda>1$ and any continuous function $\varphi(E)=\Phi'(E)$ on $\Si_{\alpha,\lambda}$,
    %the union spectrum 
    %$\Sigma_{\alpha,\lambda}^+$ 
    one has 
    \begin{equation}\label{eq:int-phi-super}
        \int \varphi(E) \, d\mmu_{\alpha,\lambda}(E) = - \left[
         \tr \left(\Phi(\Hfp_{q,\ta}) R_{\ta}\right)- \tr \left(\Phi(\Hfa_{q,0}) R_0\right)  \right].      
    \end{equation}
\end{prop}

\begin{figure}
    \centering
    \includegraphics[width=0.7\linewidth]{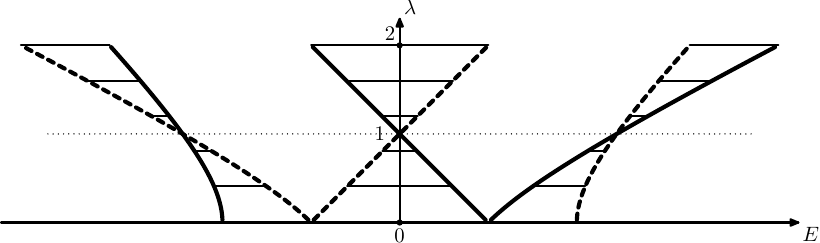} \\
    \includegraphics[width=0.7\linewidth]{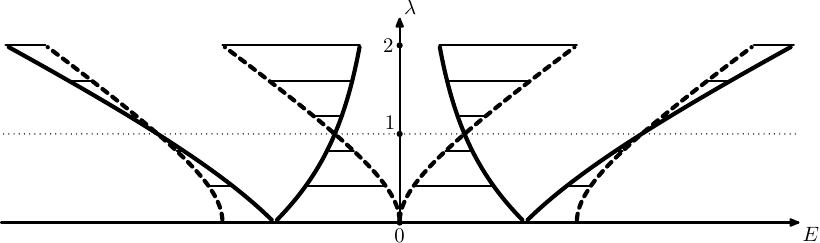}
    \caption{The set $\Si_{\frac{p}{q},\lambda}$ as a function of $\lambda$ for $\frac{p}{q}=\frac{1}{3}$ (top) %left 
    and $\frac{p}{q}=\frac{1}{4}$ (bottom). %right 
    %Red and blue 
    Dashed and solid
    curves correspond to roots of $\Hfp_{\frac{p}{q},\lambda,\ta}$ and $\Hfa_{\frac{p}{q},\lambda,0}$ respectively; horizontal lines are $q$ intervals that form~$\Si_{\frac{p}{q},\lambda}$ for the corresponding value of~$\lambda$.}
    \label{fig:graphs}
\end{figure}

\begin{remark}
    One can interpret Proposition~\ref{p:Phi-tr-super} by saying that for $\lambda>1$ the eigenvalues $\pE_{j,\ta}$ and $\aE_{j,0}$ are still endpoints of the intervals $J_j$ of $\Si_{\frac{p}{q},\lambda}$, but with the interchanged (comparing to $\lambda<1$) order; see Fig.~\ref{fig:graphs}. Such a statement also follows from our proof of the proposition.
\end{remark}

\begin{proof}[Proof of Proposition~\ref{p:Phi-tr-super}]
    Recall that $\Si_{\alpha,\lambda}= \lambda \cdot \Si_{\alpha,\frac{1}{\lambda}}$. Define the rescaled functions on $\Si_{\alpha,\frac{1}{\lambda}}$:
    \[
    \Phi_{\lambda}:=\Phi(\lambda E), \quad \varphi_{\lambda}(E)=\lambda \varphi(\lambda E) = \Phi_{\lambda}'(E).
    \]
    By definition, the left hand side of~\eqref{eq:int-phi-super} is equal to 
    \begin{multline}\label{eq:Phi-via-subcritical}
        \int_{\Si_{\alpha,\lambda}} \varphi(E) \, d\mmu_{\alpha,\lambda}(E) = \int_{\Si_{\alpha,\frac{1}{\lambda}}} \lambda \varphi(\lambda E) \, d\mmu_{\alpha,\frac{1}{\lambda}}(E) = 
        \int (\Phi_{\lambda})'(E) d\mmu_{\alpha,\frac{1}{\lambda}}
         \\           
%        =\sum_{j=1}^q (-1)^{j+1} \left(\Phi(\pE_{j,\ta; \alpha,\frac{1}{\lambda}}) - \Phi(\aE_{j,0; \alpha,\frac{1}{\lambda}})\right)         
        = \tr \left(\Phi_{\lambda}(\Hfp_{q,\ta;\alpha,\frac{1}{\lambda}}) R_{\ta}\right)
        - \tr \left(\Phi_{\lambda}(\Hfa_{q,0;\alpha,\frac{1}{\lambda}}) R_0\right)
        \\
        = \tr \left(\Phi(\lambda \Hfp_{q,\ta;\alpha,\frac{1}{\lambda}}) R_{\ta}\right) - 
        \tr \left(\Phi(\lambda \Hfa_{q,0;\alpha,\frac{1}{\lambda}}) R_0\right),
    \end{multline}
    where we have denoted $\Hfp_{q,\theta; \alpha,\lambda'}:=H_{\alpha,\lambda',\theta}|_{\mV_q^{\per}}$, 
    $\Hfa_{q,\theta; \alpha,\lambda'}:=H_{\alpha,\lambda',\theta}|_{\mV_q^{\anti}}$
    and have applied Proposition~\ref{p:Phi-tr} to pass to the second line. 

    Now, in the same way as with the classical Aubry--Andr\'e duality, we will apply the Fourier transform and check that it transforms the right hand side of~\eqref{eq:Phi-via-subcritical} into the right hand side of~\eqref{eq:int-phi-super}, conjugating $\lambda \Hfa_{q,0;\alpha,\frac{1}{\lambda}}$ and $\Hfp_{q,\ta;\alpha,\lambda}$, etc. Once performed, this (concluded by Lemma~\ref{l:Fourier} and Corollary~\ref{c:conj-lambda} below) will complete the proof.

%\todo[inline]{Say that $\lambda \Hfa_{\alpha,\frac{1}{\lambda},0}$ and $\lambda \Hfp_{\alpha,\frac{1}{\lambda},\ta}$ are conjugated to $\pm \Hfa_{\alpha,\lambda,0}$ and $\pm \Hfp_{\alpha,\lambda,\ta}$ by the Fourier transform; give reference to Fig.~\ref{fig:graphs}. VK: perhaps this will simply interchange the phases?}

    %\todo[inline]{Write the computations explicitly}

Indeed, let $(e_j)$ be the elements of the standard basis of the space $\Vp_{q}$ of $q$-periodic sequences,
\begin{equation}\label{eq:e-base}
(e_j)_n = \begin{cases}
    1, & j\equiv n \! \mod q, \\
    0, & \text{otherwise},
\end{cases}    
\end{equation}
indexed by a residue $j\in \Z/q\Z$. Take $(f_j)$ to be the Fourier base of the space $\Va_q$ of $q$-antiperiodic functions: 
\[
(f_j)_n = \exp (2\pi i \cdot (j\alpha +\ta) \cdot n);
\]
note that $q(j\alpha+\ta) \in \frac{1}{2} + \Z$ due to the choice~\eqref{eq:ta-choice} of~$\ta$, so these functions are indeed antiperiodic. Denote
%\todo{Should we use different notations for periodic/antiperiodic?} 
by $\mT_r$ the translation by $r$,
\[
(\mT_r \psi)_n = \psi_{n-r},
\]
by $\Delta_q:=\mT_1+\mT_{-1}$ the corresponding Laplacian, by $\Delta_q^{\per}$ and $\Delta_q^{\anti}$ its restrictions to~$\mV_q^{\per}$ and $\mV_q^{\anti}$ respectively, and by $D_{q,\theta}$ the diagonal operator
\[
(D_{q,\theta} \psi)_n = 2\cos 2\pi(\theta+n\alpha) \cdot \psi_n.
\]
Finally, let $\mF_q:\Vp_q \to \Va_q$ be the ``Fourier transform'' map that sends the first base into the second one,
\[
\forall j \in \Z/q\Z \quad \mF_q (e_j) = f_j.
\]
We then have the following lemma.
\begin{lemma}\label{l:Fourier}
One has
    \begin{enumerate}
        \item  $\mF_q R_{\ta} = R_0 \mF_q$,
        %$\mF_q R_0 = R_{\ta} \mF_q$
        \item $\mF_q \Delta_q^{\per} = D_{q,0} \mF_q$,
        \item $\mF_q D_{q,\ta} = \Delta_q^{\anti} \mF_q$.
    \end{enumerate}
\end{lemma}
Since 
    \[
        \lambda \Hfp_{q,\ta;\alpha,\lambda} = \lambda \Delta^{\per}_q + D_{q,\ta},
        \qquad 
        \lambda \Hfa_{q,0;\alpha,\lambda} = \lambda \Delta^{\anti}_q + D_{q,0},
    \]
once established, Lemma~\ref{l:Fourier} implies the following.
\begin{coro}\label{c:conj-lambda}
One has
\[
    \mF_q (\lambda\Hfp_{q,\ta;\alpha,\frac{1}{\lambda}}) = \Hfa_{q,0;\alpha,\lambda} \mF_q,
\]
\[
    \lambda \Hfa_{q,0;\alpha,\frac{1}{\lambda}} \mF_q= \mF_q \Hfp_{q,\ta},
\]
and hence for any continuous function $\Phi$ on $\Si_{\alpha,\lambda}$
\[
\mF_q \, \Phi(\lambda\Hfp_{q,\ta;\alpha,\frac{1}{\lambda}}) R_{\ta} \, \mF_q^{-1} = \Phi(\Hfa_{q,0;\alpha,\lambda}) R_0,
\]
\[
\mF_q^{-1}\, \Phi(\lambda\Hfa_{q,0;\alpha,\frac{1}{\lambda}}) R_{0} \, \mF_q = \Phi(\Hfp_{q,\ta;\alpha,\lambda}) R_{\ta}.
\]
\end{coro}

%Namely, fix $\theta, \theta'$ and consider the $q$-dimensional vector 

\begin{proof}[Proof of Lemma~\ref{l:Fourier}]
    \begin{multline*}        
        (\mF_q R_{\ta} e_j)_n = (\mF_q e_{-1-j})_n = \exp(2\pi i \cdot ((-1-j)\alpha+\ta) \cdot n)
        \\
        = \exp(2\pi i \cdot (-j\alpha -(\alpha-\ta))\cdot n) = 
        (\mF_q e_{j})_{-n} = (R_0 \mF_q e_j)_n,
    \end{multline*}
    where we have used that  $\alpha-\ta \equiv \ta \!\mod 1$. Next,
    \[
        (\mF_q \mT_1 e_j)_n = \exp(2\pi i \cdot ((j-1)\alpha+\ta) \cdot n) = \exp(-2\pi i \alpha n) \cdot (\mF_q e_j)_n;
    \]
    adding with the similar expression for $\mF_q \mT_{-1} e_j$, we get 
    \[
        (\mF_q \Delta_q^{\per} e_j)_n = 2\cos (2\pi \alpha n) (\mF_q e_j)_n = (D_{q,0} \mF_q e_j)_n.
    \]

    Finally,  
    \begin{multline*}
        (\mT_{1} \mF_q e_j)_n = (\mF_q e_j)_{n-1} = \exp(2\pi i \cdot (j\alpha+\ta)\cdot (n-1))
        \\
        %= \exp(2\pi i \cdot (j\alpha-\ta)) \cdot \exp(2\pi i \cdot (j\alpha-\ta)\cdot n) 
        = \exp(-2\pi i \cdot (j\alpha+\ta)) \cdot (\mF_q e_j)_n.
    \end{multline*}
    Adding with a similar equality for $\mT_{-1}\mF_q$, we get
%    \begin{multline*}
    \[
        (\Delta_q^{\anti} \mF_q e_j)_n = 2\cos (2\pi  (j\alpha+\ta)) \cdot (\mF_q e_j)_n = (\mF_q D_{q,\ta} e_j)_n.   
    \]
%    \end{multline*}
\end{proof}

As Lemma~\ref{l:Fourier} is proven, so is Corollary~\ref{c:conj-lambda}, and hence we have the equality between the traces of $\mF_q$-conjugated operators:
\[
\tr\left( \Phi(\lambda\Hfp_{q,\ta;\alpha,\frac{1}{\lambda}}) R_{\ta} \right) = \tr\left( \Phi(\Hfa_{q,0;\alpha,\lambda}) R_0 \right),
\]
\begin{equation}\label{eq:0-to-ta}
\tr\left( \Phi(\lambda\Hfa_{q,0;\alpha,\frac{1}{\lambda}}) R_{0} \right) = \tr\left( \Phi(\Hfp_{q,\ta;\alpha,\lambda}) R_{\ta}\right).    
\end{equation}
Joining it with~\eqref{eq:Phi-via-subcritical} gives us the desired~\eqref{eq:int-phi-super}, and thus completes the proof of Proposition~\ref{p:Phi-tr-super}.
\end{proof}

\section{Infinite dimensional operators}\label{s:infinite}
\subsection{Moments $c_{2k}(\alpha,\lambda)$}\label{s:c-moments}

Note that Propositions~\ref{p:Phi-tr} and~\ref{p:Phi-tr-super} already allow to establish Theorem~\ref{t:poly} for rational values of $\alpha$ with sufficiently large (compared to~$k$) denominators. 

\begin{figure}
    \centering
    \includegraphics[width=0.5\linewidth]{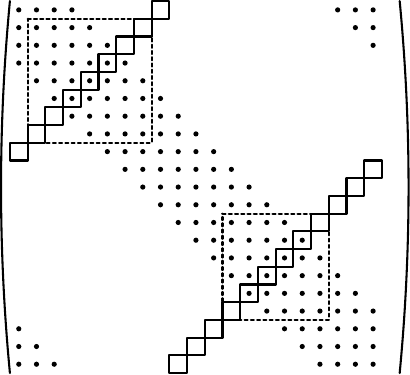}
    \caption{Matrix of $(H^{\mathrm{per}/\mathrm{anti}}_{\frac{p}{q},\lambda,\theta})^{2k+1}$ for $q\gg k$. Nonzero elements are shown by dots, anti-diagonal elements marked by squares; coordinates are cyclically shifted to put both intersections between the diagonal and the anti-diagonal away from the edges of the matrix. Areas around these intersections are shown by dashed squares.}
    \label{fig:intersections}
\end{figure}

 Namely, the moment $c_{2k}(\alpha,\lambda)$ corresponds (in the sense of Proposition~\ref{p:Phi-tr}) to $\Phi(E)=\frac{E^{2k+1}}{2k+1}$. The matrices, whose traces are being calculated, are (up to the factor $\frac{1}{2k+1}$) of the form $H^{2k+1} R$ for one of $H=\Hfa_{q,0},\Hfp_{q,\ta}$ and the corresponding reflection $R=R_0,R_{\ta}$. The traces can be thus seen, up to the signs coming from the (anti)periodicity, as the sum of the \emph{anti}-diagonal elements of $H^{2k+1}$. Now, the non-zero elements of these matrices are located in $2k+1$-neighbourhood of the main diagonal, and anti-diagonal intersects this region in \emph{two} places (that are separated by about $q/2$ indices from each other); see Fig.~\ref{fig:intersections}. The behaviour of the finite dimensional matrices $H^{2k+1}$ around these intersection points can then be studied separately, and described in terms of an infinite matrix of the operator $H$ with the corresponding phase, acting on $\ell^2(\Z).$

 To obtain such a description, first let us consider the set 
\begin{equation}\label{eq:Th-def}
\Th_{\alpha}:=\left\{0, \frac{\alpha}{2}, \frac{1}{2}, \frac{1+\alpha}{2}\right\}    
\end{equation} 
of phases; note that each phase $\theta\in \Th_{\alpha}$ one has $-\theta\equiv \theta \mod (\Z + \alpha \Z)$. These four points correspond to the regions of intersection between the diagonal and the anti-diagonal that one can associate to the matrices in the right hand side of~\eqref{p:Phi-tr}; there are four such regions, two for each of the traces in the right hand side. For each of the phases $\theta\in \Th_{\alpha}$, we define the symmetry on $\ell^2(\Z)$ by
\begin{equation}\label{eq:R-12}
R_0, R_{\frac{1}{2}}:(\psi_n)\mapsto (\psi_{-n}), \quad 
R_{\fa},R_{\frac{1+\alpha}{2}}:(\psi_n)\mapsto (\psi_{-1-n});
\end{equation}
these symmetries then commute with the corresponding $H_{\alpha,\lambda,\theta}$, and this definition also agrees  with~\eqref{eq:R0} and~\eqref{eq:Rfa}.

We then have the following proposition.

\begin{prop}\label{p:moment-trace}
    Assume that $\alpha=\frac{p}{q}$, where $q>4k+2$, and $\lambda<1$. Then
    \begin{multline}\label{eq:moment-traces}
        (2k+1) c_{2k}(\alpha,\lambda) = \tr H_{\alpha,\lambda,\frac{1+\alpha}{2}}^{2k+1} R_{\frac{1+\alpha}{2}} + \tr H_{\alpha,\lambda,\fa}^{2k+1} R_{\fa} +
        \\
        +\tr H_{\alpha,\lambda,\frac{1}{2}}^{2k+1} R_{\frac{1}{2}}- \tr H_{\alpha,\lambda,0}^{2k+1} R_0,
    \end{multline}
    where the traces in the right hand side are understood as the sum of the diagonal elements in the standard base of $\ell^2(\Z)$. For $\lambda>1$, the right hand side of~\eqref{eq:moment-traces} changes the sign.
\end{prop}

\begin{remark}
    The traces in right hand side of~\eqref{eq:moment-traces}  again can be interpreted as sums of anti-diagonal elements of the corresponding powers of $H_{\alpha,\lambda,\theta}$ (see~\eqref{eq:H-R-0},~\eqref{eq:H-R-12}, where the placements of these elements are highlighted). From this interpretation, it is easy to see that only a finite number of elements in the corresponding sums are non-zero.
\end{remark}

\begin{equation}\label{eq:H-R-0}
H_{\alpha,\lambda,0}=\begin{pmatrix}
\ddots & 1 & 0 & \phantom{\ddots} & \framebox{$0$} \\
1 & V_{-\alpha} & 1 & \framebox{$0$} & \phantom{\ddots}& \\
0  & 1 & \framebox{$V_0$} & 1 & 0  \\
\phantom{\ddots} & \framebox{$0$}  & 1 & V_{\alpha} & 1 \\
\framebox{$0$} & \phantom{\ddots}  & 0  & 1 &  \ddots 
\end{pmatrix},
\quad  H_{\alpha,\lambda,\frac{1+\alpha}{2}}=
\begin{pmatrix}
  \ddots & 1 & $0$ & \framebox{$0$}& \\
 1 & V_{\frac{1-\alpha}{2}} & \framebox{$1$} & 0   \\
0  & \framebox{$1$} & V_{{\frac{1+\alpha}{2}}} & 1  \\
 \framebox{$0$}  & 0  & 1 & \ddots 
\end{pmatrix},
\end{equation}
\begin{equation}  \label{eq:H-R-12}  
H_{\alpha,\lambda,\frac{1}{2}}=\begin{pmatrix}
\ddots & 1 & 0 & \phantom{\ddots} & \framebox{$0$} \\
1 & V_{\frac{1}{2}-\alpha} & 1 & \framebox{$0$} & \phantom{\ddots}& \\
0  & 1 & \framebox{$V_{\frac{1}{2}}$} & 1 & 0  \\
\phantom{\ddots} & \framebox{$0$}  & 1 & V_{\frac{1}{2}+\alpha} & 1 \\
\framebox{$0$} & \phantom{\ddots}  & 0  & 1 &  \ddots 
\end{pmatrix},
\quad  H_{\alpha,\lambda,\frac{\alpha}{2}}=
\begin{pmatrix}
  \ddots & 1 & $0$ & \framebox{$0$}& \\
 1 & V_{-\frac{\alpha}{2}} & \framebox{$1$} & 0   \\
0  & \framebox{$1$} & V_{\frac{\alpha}{2}} & 1  \\
 \framebox{$0$}  & 0  & 1 & \ddots 
\end{pmatrix},
\end{equation}

%\todo[inline]{Do we want to compare the four-traces formulas with~\eqref{eq:AvMS-Leb}? Maybe not, as these are different }

 Proposition~\ref{p:moment-trace} is a particular case of Theorem~\ref{t:four-traces} below. Namely, it turns out that the conclusion of Proposition~\ref{p:moment-trace} can be generalized to the integral of an arbitrary $\varphi(E)=\Phi'(E)$, where the function $\Phi$ is analytic in some neighbourhood of $\Spl_{\alpha,\lambda}$. Thus, instead of a complete proof of this proposition, we present only a sketch that clarifies the approach, and formally establish it with the proof of Theorem~\ref{t:four-traces}.

\begin{proof}[Sketch of the proof of Proposition~\ref{p:moment-trace}]
    Consider first the trace $\tr \left( (\Hfp_{q,\ta})^{2k+1} R_{\ta}\right)$ in the right hand side of~\eqref{eq:Phi-trace}. Considering the coordinates as residues modulo $q$, one has
    \begin{equation}\label{eq:sum-Z-q}
        \tr \left( (\Hfp_{q,\ta})^{2k+1} R_{\ta}\right)=
        \sum_{i\in \Z/q\Z} (\Hfp_{q,\ta})^{2k+1}_{i,-1-i}.
    \end{equation}
    Now, the element $(\Hfp_{q,\ta})^{2k+1}_{i,j}$ can be nonzero only if the distance between $i$ and $j$ in $\Z/q\Z$ is at most $2k+1$. Hence, nonzero summands in the right hand side of~\eqref{eq:sum-Z-q} correspond either to $|i|\le k+1$ or to $|i-\frac{q}{2}|\le k+1$. In the former case we get exactly the same summands as for $\tr \left( (H_{\frac{p}{q},\lambda,\ta})^{2k+1} R_{\ta}\right)$, in the latter the same as for $\tr \left( (H_{\frac{p}{q},\lambda,\theta'})^{2k+1} R_{\theta'}\right)$, where $\theta'\in \Th_{\alpha}$ is another point in the same orbit of rotation by $\alpha$ as~$\ta$.

    In the same way, one handles the trace 
    $\tr \left( (\Hfa_{q,0})^{2k+1} R_{0}\right)$, getting two other traces that appear in the right hand side of~\eqref{p:moment-trace}. The trace $\tr \left( (H_{\frac{p}{q},\lambda,0})^{2k+1} R_{0}\right)$ appears with the same minus sign as $\tr \left( (\Hfa_{q,0})^{2k+1} R_{0}\right)$ in the right hand side of~\eqref{eq:Phi-trace}, while the fourth trace appears with the plus sign due to the antiperiodicity that changes the sign of its contribution.
\end{proof}

As we will see below, actually the statement of Proposition~\ref{p:moment-trace} holds without the assumption $q>4k+2$. And this is where polynomial formulas for the moments, declared in Theorem~\ref{t:poly}, are coming from: all the four traces in the right hand side of~\eqref{eq:moment-traces} are polynomial in $\lambda$ with coefficients that are polynomials of $\cos \pi \alpha$, that with some additional consideration can be seen to be actually polynomials of~$\cos 2\pi \alpha$. One can \emph{compute} the polynomials for the second and fourth moment using the right hand side of~\eqref{eq:moment-traces}, though Proposition~\ref{p:moment-trace} provides the proof of this answer only for rational values of~$\alpha$; it is Theorem~\ref{t:four-traces} below that will establish the correctness of this answer in full generality. Also, the following lemma allows to simplify the computations.

% \todo[inline]{Examples discussion: lemma \& corollary below already allow to guess the polynomials for the moments.}
\begin{lemma}\label{p:simplifying}
    For any $\alpha=\frac{p}{q}$ and $\lambda\in\R$ one has 
    \begin{equation}\label{eq:equal}        
        \tr H_{\alpha,\lambda,\frac{1}{2}}^{2k+1} R_{\frac{1}{2}} =- \tr H_{\alpha,\lambda,0}^{2k+1} R_0,
        \qquad
        \tr H_{\alpha,\lambda,\frac{\alpha}{2}}^{2k+1} R_{\frac{\alpha}{2}} = \tr H_{\alpha,\lambda,\frac{1+\alpha}{2}}^{2k+1} R_{\frac{1+\alpha}{2}}.
    \end{equation}
\end{lemma}
\begin{proof}
    As the potential satisfies $V(\theta+\frac{1}{2})=-V(\theta)$, for every $\theta$ the diagonal elements of the matrices $H_{\alpha,\lambda,\theta}$ and  $H_{\alpha,\lambda,\theta+\frac{1}{2}}$ differ by the sign, while the off-diagonal elements  coincide. Hence, for any $\theta$, $k$ and any indices $i,j\in \Z$ one has 
    \[
        (H_{\alpha,\lambda,\theta}^{2k+1})_{i,j}= (-1)^{i-j+1} \cdot (H_{\alpha,\lambda,\theta+\frac{1}{2}}^{2k+1})_{i,j}.
    \]
    Now it remains to notice that for the first pair of matrices in~\eqref{eq:equal} we are summing elements with coordinates $i$ and $j$ of the same parity (hence the difference of the sign), while for the second pair the coordinates are of opposite parity (hence the corresponding elements coincide).
\end{proof}

\begin{figure}
    \centering
    \includegraphics[width=0.35\linewidth]{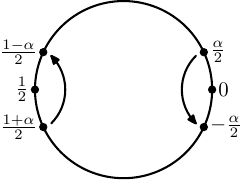}
    \caption{Four $\mR$-invariant orbits of rotation by~$\alpha$.}
    \label{fig:circle-alpha}
\end{figure}

%  {Roughly speaking, in the standard base of the space $\mV_q^{\per}$ the trace of the product of an operator with the symmetry $R_{\ta}$ is the sum of its \emph{anti}-diagonal elements. The matrix of $H_{\alpha,\lambda,\ta}$ is three-diagonal, hence its $2k+1$-st power is concentrated on the diagonals within range $2k+1$ of the main one. However, the diagonal and the anti-diagonal on the torus intersect in two points. Hence, the trace $\tr \left(\Phi(\Hfa_{q,0}) R_0\right)$ will be equal to the sum of two sums of anti-diagonal elements of infinite-dimensional matrices, $H_{\alpha,\lambda,\ta}^{2k+1}$ and $H_{\alpha,\lambda,\theta'}^{2k+1}$, where $\theta'$ is another point of the same orbit. 
% }

%\vspace{1cm}

% Now, for each of the values $\theta\in \{0, \frac{\alpha}{2}, \frac{1}{2}, \frac{1+\alpha}{2}\}$, consider the corresponding action of $\mR$ on $\ell^2(\mO_{\theta})$; identifying these spaces with $\ell^2(\Z)$, we get the operators that commute with the corresponding Hamiltonians $H_{\alpha,\lambda,\theta}$. Namely, 

\begin{example}\label{ex:compute-second}
    Let us find the expression for the second moment $c_2(\alpha,\lambda)$, using rational $\alpha=\frac{p}{q}$ for sufficiently large $q$ for Proposition~\ref{p:moment-trace} to be applicable. The $3\times 3$ submatrix with indices $i=-1,0,1$ of $H_{\alpha,\lambda,0}$ equals
    \[
        A=\begin{pmatrix}
        2\lambda\cos 2\pi\alpha & 1 & \framebox{$0$} \\
        1 & \framebox{$2\lambda$} & 1  \\
        \framebox{$0$}  & 1 & 2\lambda\cos 2\pi\alpha  \\
        \end{pmatrix};
    \]
    no other elements contribute to $\tr H_{\alpha,\lambda,0}^3 R_0$. Now, the elements of $A^3$ that contribute to this trace are
    \[
        (A^3)_{1,-1}=(A^3)_{-1,1}=2\lambda\cos 2\pi\alpha + 2\lambda + 2\lambda\cos 2\pi\alpha = 2\lambda (1+2\cos 2\pi\alpha),
    \]
    and 
    \[
        (A^3)_{0,0}=(2\lambda)^3 + 4\cdot (2\lambda\cdot 1 \cdot 1) + 2\cdot (2\lambda\cos 2\pi\alpha \cdot 1\cdot 1) = 8\lambda^3 + 2\lambda\cdot (4+2\cos 2\pi\alpha),
    \]
    thus 
    \begin{equation}\label{eq:tr-3-odd}
        \tr H_{\alpha,\lambda,0}^3 R_0 = (A^3)_{-1,1} + (A^3)_{0,0} + (A^3)_{1,-1} = 8\lambda^3 + 12\lambda (1+\cos 2\pi\alpha).        
    \end{equation}
    This trace comes to the right hand side of~\eqref{eq:moment-traces} with a minus sign, and due to Lemma~\ref{p:simplifying}, one gets an equal contribution from $\tr H_{\alpha,\lambda,\frac{1}{2}}^3 R_{\frac{1}{2}}$. Next, one can perform a similar computation using the $4\times 4$ submatrix of $H_{\alpha,\lambda,\frac{\alpha}{2}}$, corresponding to the indices $i=-2,-1,0,1$,
    \[
    A'=\begin{pmatrix}
        2\lambda \cos 3\pi \alpha & 1 & 0 & \framebox{$0$}& \\
        1 & 2\lambda \cos \pi \alpha & \framebox{$1$} & 0   \\
        0  & \framebox{$1$} & 2\lambda \cos \pi \alpha & 1  \\
        \framebox{$0$}  & 0  & 1 & 2\lambda \cos 3\pi \alpha 
    \end{pmatrix},
    \]
    obtaining 
    \[
        (A')^3_{-2,1}=(A')^3_{1,-2}=1,
    \]
    \[
        (A')^3_{-1,0}=(A')^3_{0,-1}=3+ 3\cdot 4\lambda^2 \cos^2 \pi \alpha
    \]
    and thus
    \begin{equation}\label{eq:tr-3-even}
        \tr H_{\alpha,\lambda,\frac{\alpha}{2}}^3 R_{\frac{\alpha}{2}}=\tr H_{\alpha,\lambda,\frac{1+\alpha}{2}}^3 R_{\frac{1+\alpha}{2}}=8 + 12\lambda^2\cdot (1+\cos 2\pi\alpha).
    \end{equation}
    One thus finally gets 
    \[
        \frac{3}{16} c_{2}(\alpha,\lambda) = 1-\frac{3}{2}\lambda (1+\cos 2\pi\alpha) + \frac{3}{2}\lambda^2 (1+\cos 2\pi\alpha) -\lambda^3,
    \]
    as declared in Example~\ref{ex:2-4}. Finally, note that yet another manifestation of Aubry-Andr\'e duality allows to simplify the computation even further: as the reader might have already noticed from comparing~\eqref{eq:tr-3-odd} and~\eqref{eq:tr-3-even}), one has
    \begin{equation}\label{eq:A-A-duality-bis}
        \tr H_{\alpha,\lambda,\frac{\alpha}{2}}^{2k+1} R_{\frac{\alpha}{2}} = - \lambda^{2k+1} \tr H_{\alpha,\frac{1}{\lambda},0}^{2k+1} R_{0}.        
    \end{equation}
    Indeed, note first that it suffices to verify~\eqref{eq:A-A-duality-bis} for the case of a sufficiently large even denominator $q$: this is an identity between polynomials in $\lambda$ and $\cos 2\pi\alpha$, and it suffices to check it on a dense set of arguments.
    For such $q$, in terms of the proof of Proposition~\ref{p:moment-trace} one has $\{\ta,\theta'\}=\{\frac{\alpha}{2},\frac{1+\alpha}{2}\}$ and $\theta''=\frac{1}{2}$, and thus from Lemma~\ref{p:simplifying}
    \[
    \tr \left( (\Hfp_{q,\ta})^{2k+1} R_{\ta}\right) = \tr \left( H_{\alpha,\lambda,\ta}^{2k+1} R_{\ta}\right)+ \tr \left( H_{\alpha,\lambda,\theta'}^{2k+1} R_{\theta'}\right) 
    = 2\tr \left( H_{\alpha,\lambda,\frac{\alpha}{2}}^{2k+1} R_{\frac{\alpha}{2}}\right),
    \]
    \[
    \tr \left( (\Hfa_{q,0;\alpha,\frac{1}{\lambda}})^{2k+1} R_{0} \right) = \tr \left( H_{\alpha,\frac{1}{\lambda},0}^{2k+1} R_{0}\right)- \tr \left( H_{\alpha,\frac{1}{\lambda},\theta''}^{2k+1} R_{\theta''}\right) 
    = 2\tr \left( H_{\alpha,\frac{1}{\lambda},0}^{2k+1} R_{0}\right).
    \]
    Now, from Corollary~\ref{c:conj-lambda} and Eq.~\eqref{eq:0-to-ta} we have
    \begin{multline*}
        - \lambda^{2k+1} \tr \left( (\Hfa_{q,0;\alpha,\frac{1}{\lambda}})^{2k+1} R_{0} \right) = 
        \tr \left( (\lambda \Hfa_{q,0;\alpha,\frac{1}{\lambda}})^{2k+1} R_{0} \right) 
        \\ 
        = 
        \tr\left( (\Hfp_{q,\ta;\alpha,\lambda})^{2k+1}R_{\ta}\right),    
    \end{multline*}
    and the equality~\eqref{eq:A-A-duality-bis} follows.
    
    % $\tr H_{\alpha,\lambda,0}^{2k+1} R_0$ as a function of $\lambda$ is odd (any product of $2k+1$ matrix elements contributing to this trace has an odd number of elements on the diagonal), while for the same reason $\tr H_{\alpha,\lambda,\frac{\alpha}{2}}^{2k+1} R_{\frac{\alpha}{2}}$ is even. Joining it with Hence, 
    % %.by using Andr\'e-Aubry duality from Corollary~\ref{c:conj-lambda}: one has 
    % \[
    %     \tr H_{\alpha,\lambda,\ta}^{2k+1} R_{\ta} = - \lambda^{2k+1} \tr H_{\alpha,\frac{1}{\lambda},0}^{2k+1} R_{0}.
    % \]
    % (The reader might have already noticed a manifestation of this duality by comparing~\eqref{eq:tr-3-odd} and~\eqref{eq:tr-3-even}.)
    % by noticing that $\tr H_{\alpha,\lambda,0}^{2k+1} R_0$ as a function of $\lambda$ is odd (any product of $2k+1$ matrix elements contributing to this trace has an odd number of elements on the diagonal), while $\tr H_{\alpha,\lambda,\frac{\alpha}{2}}^{2k+1} R_{\frac{\alpha}{2}}$ is even. Thus, the contribution of $\tr H_{\alpha,\lambda,\frac{\alpha}{2}}^{2k+1} R_{\frac{\alpha}{2}}$ can be reconstructed using Andr\'e-Aubry duality 
    % %
\end{example}

%\todo[inline]{Phrase here that actually the formula works for all $\lambda$, and this is a theorem below, but to state it, we need to check that RHS is well-defined}

\begin{remark}
    Consider the action of the reflection map $\mR:\theta\mapsto -\theta$ on the circle of phases (see Fig~\ref{fig:circle-alpha}); it preserves the potential $V(\theta)=2\lambda \cos 2\pi \theta$. Consider then its action on the set of orbits of rotation by~$\alpha$.
    For rational values of $\alpha=\frac{p}{q}$, considered above, there are two orbits (the one of~$0$ and of~$\ta$), that are preserved by~$\mR$, and each of them corresponds to two trace terms in~\eqref{eq:moment-traces}. But when $\alpha$ becomes irrational (as it will be in Theorem~\ref{t:four-traces} below), we get exactly four different orbits fixed by~$\mR$, exactly those passing through the points of~$\Th_{\alpha}$, and each these orbits corresponds to one of the trace terms in the right hand side of~\eqref{eq:four-traces} below.
\end{remark}

% \emb{Consider the action of the reflection map $\mR:\theta\mapsto -\theta$ on the circle of phases; it preserves the potential $V(\theta)=2\lambda \cos 2\pi \theta$. For an \emph{irrational} $\alpha$, among the orbits 
% \[
% \mO_{\theta}=\{\theta + \alpha n \mid n\in \Z \} \subset \R/2\pi\Z
% \]
% of the rotation by $\alpha$ there are four ones that are sent by $\mR$ to itself (see Fig.~\ref{fig:circle-alpha}). Indeed, solving the congruence 
% \[
% \theta\equiv -\theta \mod \alpha \Z + \Z, 
% \]
% leads to the orbits of phases
% \[
% \theta\in \left\{0, \frac{\alpha}{2}, \frac{1}{2}, \frac{1+\alpha}{2}\right\}.
% \]}

 \begin{remark}
    Formula~\eqref{eq:moment-traces} (as well as its counterpart~\eqref{eq:four-traces} in Theorem~\ref{t:four-traces} below) looks miraculous: why exactly these ``traces'' of these (infinite) matrices appear? Can it be considered as one of a huge variety of traces formulae, for instance, as a distant relative of Lefshets formula, or of any kind of index theorems? It would be interesting to have even a hand-waving explanation; for the moment the authors have none.
\end{remark}

 \subsection{Decay estimates and analyticity of elements of~$\Phi(H)$}\label{s:decrease}
 We will need the following proposition:
 \begin{prop}\label{p:Phi-bounds}
     Let function $\Phi$ be analytic in some neighbourhood of $\Sigma^+_{\alpha,\lambda}$. Then for every phase $\theta$, the elements of $\Phi(H_{\alpha,\lambda,\theta})$ decay exponentially with the distance from the diagonal: there exists $C_0, c_0>0$:
     \[
        \forall j,l \quad \left|(\Phi(H_{\alpha,\lambda,\theta}))_{jl}\right| \le C_0 \exp (-c_0 |j-l|).
     \]
     Moreover, this estimate is uniform in $\theta$ and in some neighbourhood of the initial point $(\alpha,\lambda)$. Finally, all these elements depend analytically on $(\lambda,\theta)$ and $C^{\infty}$-smoothly on $\alpha$, and all their derivatives satisfy the exponential decay estimates (possibly each with different constants~$c_0,C_0$).
 \end{prop}

We will deduce Proposition~\ref{p:Phi-bounds} via Cauchy's formula from the following statement, that could be seen as a variation of (one-dimensional) version of Combes--Thomas estimates (\cite{CT}; also, see, e.g. \cite[{Theorem~11.2}]{K} for a modern exposition). More specifically, the next lemma states the same conclusions for a particular case of $\Phi_z(E)=1/(E-z)$, corresponding to the Green's function (or resolvent)
\begin{equation}\label{eq:Green-def}
    G_{\alpha,\lambda,\theta}(z)=(H_{\alpha,\lambda,\theta}-z\cdot \Id)^{-1}.
\end{equation}

\begin{lemma}\label{l:Green}
    For any $\alpha,\lambda$ and any $E_0 \notin \Spl_{\alpha,\lambda}$, for some $C_1,c_1>0$ the elements of the Green's function satisfy uniform estimate
    \begin{equation}
        \forall j,l\in \Z, \quad \forall \theta\in \Sc \quad 
        \left|(G_{\alpha,\lambda,\theta}(E_0))_{jl}\right| \le C_1 \exp (-c_1 |j-l|);        
    \end{equation}
    moreover, the constants $C_1,c_1$ can be chosen uniformly in some neighbourhood of $(\alpha,\lambda,E_0)$ in $\R\times\C\times \C$. Finally, the elements of the Green's function depend analytically on $\lambda,E_0,\theta$ and $C^{\infty}$-smoothly on~$\alpha$, % and~$\theta$, 
    and for any derivative a similar uniform exponential decay also holds in some neighbourhood of $(\alpha,\lambda,E_0)$, possibly with different constants.
    %, and similar estimates (with other constants) hold for their derivatives.
\end{lemma}

\begin{proof}[Proof of Lemma~\ref{l:Green}]
    Consider the skew product 
    \begin{equation}\label{eq:skew}
        F_{\lambda,E}: (\Sc)^2\times \C^2 \to (\Sc)^2\times \C^2,  \quad F_{\lambda,E}(\alpha,\theta;v) = (\alpha,\theta+\alpha; \Pi_{\theta,E} (v) ),
    \end{equation}
    where $\Pi_{\theta,E}$ are given by~\eqref{eq:Pi-def}. Then, as $E_0\notin \Spl_{\alpha,\lambda}$, due to Johnson's Theorem the skew product~\eqref{eq:skew} is a hyperbolic $\SL(2,\C)$-cocycle \cite[Theorem 5.2.8]{DF24} (also, see~\cite{Jo86} for the original paper where Johnson's Theorem was initially stated, or \cite{Zh20} for a modern exposition), that is, there exists a splitting of fibers $\C^2 = E^s_{\alpha,\theta} \oplus E^s_{\alpha,\theta}$, invariant under the dynamics, such that the vectors from $E^s$ are exponentially contracted, while the vectors from $E^u$ are exponentially expanded.

    Now, we have 
    \[
        (G_{\alpha,\lambda,\theta}(E_0))_{jl} = 
        (G_{\alpha,\lambda,\theta+l\alpha}(E_0))_{j-l,0},
    \]
    so (as the desired estimate should be uniform in $\theta$) we can assume without loss of generality that $l=0$. Take the vector $\psi=(\psi_n)_{n\in\Z}\in \ell^2(\Z)$ to be the $0$-th column of the Green's function matrix, 
    \[
    \psi_n = (G_{\alpha,\lambda,\theta}(E_0))_{n, 0};
    \]
    by definition of the Green's function, it satisfies the equation 
    \begin{equation}\label{eq:H-E0}
        H_{\alpha,\lambda,\theta}\psi - E_0 \psi = \delta_{n 0},         
    \end{equation}
    where 
    \[
    \delta_{n0}=\begin{cases} 1 & n=0 \\ 0 & n\neq 0 \end{cases}
    \]
    is the Kronecker symbol. 
    Equation~\eqref{eq:H-E0} is equivalent to the recurrent relation
    \begin{equation}\label{eq:H-E0-bis}
        \psi_{n+1} = (E_0-V(\theta+n\alpha))\psi_n - \psi_{n-1} + \delta_{n 0}.
    \end{equation}
    Now, consider the sequence of vectors $v_n\in \C^2$ defined by $v_n=\left( {\psi_{n+1} \atop \psi_n} \right)$; then,~\eqref{eq:H-E0-bis} can 
    be rewritten as %is equivalent to
    %and hence to
    \begin{equation}\label{eq:v-rec}
        v_{n}= \Pi_{E_0,\theta+(n-1)\alpha} v_{n-1}+ \delta_{n0} \left( {1 \atop 0} \right).        
    \end{equation}
    Now, take the decomposition
    \begin{equation}\label{eq:v-delta}
        \left( {1 \atop 0} \right)= v_0^s + v_0^u, \quad v_0^s\in E_{\alpha,\theta}^s, \quad v_0^u\in E_{\alpha,\theta}^u;
    \end{equation}
    The solution to~\eqref{eq:v-rec} can be then be written as
    %Then the solution to~\eqref{eq:v-rec} can be written as
    \begin{equation}\label{eq:v-n-product}
        v_n = \begin{cases}
            \Pi_{E_0,\theta+(n-1)\alpha} \dots \Pi_{E_0,\theta} v_0^s, & n \ge 0, \\
            -\Pi_{E_0,\theta-n \alpha}^{-1} \dots \Pi_{E_0,\theta-\alpha}^{-1} v_0^u, & n < 0.    
        \end{cases}        
    \end{equation}
    The exponential decay (Combes--Thomas) estimates then follow from the hyperbolicity of the cocycle, as both forward iterations of $v_0^s \in E_{\alpha,\theta}^s$ and backward iterations of $v_0^u\in E_{\alpha,\theta}^u$ decay exponentially. 

    Next, note that the stable and unstable directions can be found by the graph transform, acting on the maps from the product of (complex) neighborhoods of initial $\lambda,E_0$ and of the phase space $U^{\lambda}\times U^{E_0}\times U_{\eps}(\R/Z)$ to the set of (complex) directions $\CP^1$. Uniform hyperbolicity of the cocycle implies that (for sufficiently small neighbourhoods) this transformation and its inverse are contracting in $C^0$-neighbourhoods of the unstable and stable directions respectively. This implies that the decomposition~\eqref{eq:v-delta} depends analytically on $\lambda$ and $\theta$, see \cite[Lemma 2.1]{A} and \cite[Theorem 6.1]{AJS}.

    Moreover, for every $r$, consider the action of graph transform on the space $r$-jets in $\alpha$, holomorphic in $(\lambda,E_0,\theta)$; then, due to zero Lyapunov exponents of the transformation in the base, some its iteration will again be contracting. This implies that the stable and unstable directions depend $C^{\infty}$-smoothly on $\alpha$, with each $r$-th derivative being analytic in~$(\lambda,E_0,\theta)$. As an overview remark, this argument is in the spirit of normal hyperbolicity: the base transformation $(\alpha,\lambda,E_0,\theta)\mapsto (\alpha,\lambda,E_0,\theta+\alpha)$ has zero Lyapunov exponent, as opposed to the transverse exponential contraction for the invariant graph, thus implying existence of all the derivatives. For general real-valued transformations (see~\cite{HPS} for a detailed exposition of the theory of normal hyperbolicity), it would imply only $C^{\infty}$-smoothness in $(\alpha,\theta)$, but as the transformation is complex-analytic in $(\lambda,E_0,\theta)$ (and a complex neighbourhood of the phase circle is invariant), the dependence on $\theta$ is actually complex-differentiable, and hence analytic.
    %; for the latter, note that the dynamics can be extended to some complex neighbourhood $U_{\eps}(\R/\Z)$ of the phase circle, and 
    %this neighbourhood is preserved by the rotation by~$\alpha$. 
    %Hence, the decomposition~\eqref{eq:v-delta} also depends analytically on $\lambda$ and~$\theta$. 
    
    % At the same time, the projectivisation image of the unstable (or stable) direction forms a normally hyperbolic manifold, with transverse contraction and zero Lyapunov exponent along $(\alpha,\theta)$, and hence it is $C^{\infty}$-smooth in~$\alpha$ (see \cite{HPS} for the detailed exposition of the theory of normal hyperbolicity).
    %\todo[inline]{Add something?}

    Now, it is a standard argument that a uniformly hyperbolic $\SL(2,\C)$-cocycle is conjugate to a diagonal one. In our case, one could choose this conjugacy to be analytic in $(\lambda,E_0,\theta)$ and smooth in~$\alpha$. To prove it, one starts by showing that in a product of (perhaps, smaller) neighbourhoods of $\lambda, E_0$ and the phase circle $\R/\Z$ there exists a change of variables $B=B_{\alpha,\lambda,E_0,\theta}$ of the same regularity that sends stable and unstable directions to the horizontal $\C\times \{0\}$ and vertical $\{0\}\times \C$ directions respectively. Indeed, the images 
    \[
        \{[E_{\alpha,\theta}^s] \mid \theta\in \R/\Z \}, 
        \{[E_{\alpha,\theta}^u] \mid \theta\in \R/\Z \}\subset \CP^1
    \]
    are images of the phase circle $\R/\Z$ under smooth maps, and thus are of measure zero. Hence, there exists a coefficient $s$ such that the linear function 
    \[
    L:\C^2\to \C, \quad L(z_1,z_2)=z_1 + s z_2
    \]
    does not vanish on the spaces $E^s_{\alpha,\theta}, E^s_{\alpha,\theta}$ for any $\theta\in\R/\Z$. Fix such a function~$L$. By continuity, the same statement holds for $\alpha,\lambda,E_0$ in some neighbourhood of their initial values, as well as in some complex $\eps$-neighbourhood of the phase circle.

    Now, fix the vectors 
    \[
    v^s_{\alpha,\theta}\in E^s_{\alpha,\theta}, \quad v^u_{\alpha,\theta}\in E^u_{\alpha,\theta}
    \]
    by the relations
    \[
    L(v^s_{\alpha,\theta})= L(v^u_{\alpha,\theta})=1.
    \]
    Then $v^s_{\alpha,\theta},v^u_{\alpha,\theta}$ are nonzero vectors in the corresponding directions, and depend on the parameters as regularly as the directions do. Take these as the columns of the matrix of the coordinate change $B_{\alpha,\theta}$. Then applying $B_{\alpha,\theta}^{-1}$ sends the stable and unstable directions to the horizontal and the vertical ones. As the decomposition into these directions is invariant, in the new coordinates the cocycle becomes diagonal:
    \begin{equation}\label{eq:B-conj}        
        B_{\alpha,\theta+\alpha} \Pi_{E_0,\theta} B_{\alpha,\theta}^{-1} = 
        \begin{pmatrix}
            \kappa^s_{\alpha,\theta} & 0 \\
            0 & \kappa^u_{\alpha,\theta}
        \end{pmatrix}.
    \end{equation}
    Also, due to the hyperbolicity of the cocycle, upon further conjugating with a diagonal ($\theta$-dependent) matrix, we can assume that 
    \[
        \exists \rho<1 : \forall \theta\in\Sc \quad |\kappa^s_{\alpha,\theta}|<\rho, \quad |\kappa^u_{\alpha,\theta}|>\rho^{-1}.
    \]
    
    Now, applying $B_{\alpha,\theta}$ to the representation~\eqref{eq:v-delta}, we get 
    \[
        B_{\alpha,\theta} v^s_0 = 
        \begin{pmatrix}
            t^s_{\alpha,\theta} \\ 0
        \end{pmatrix},
        \quad 
        B_{\alpha,\theta} v^s_0 = 
        \begin{pmatrix}
            0 \\ t^u_{\alpha,\theta}
        \end{pmatrix},
    \]
    with the coordinate functions $t^s, t^u$ that are $C^{\infty}$-smooth in $\alpha$ (and analytic in $\lambda,\theta$).
    %applying~\eqref{eq:B-conj} to
    Rewriting~\eqref{eq:v-n-product} using~\eqref{eq:B-conj}, we get for $n>0$
    \begin{equation}\label{eq:v-n-change}
        v_n = B_{\alpha,\theta+n\alpha}^{-1} \begin{pmatrix}
            \kappa^s_{\alpha,\theta+(n-1)\alpha} \dots \kappa^s_{\alpha,\theta +\alpha} \kappa^s_{\alpha,\theta} t^s_{\alpha,\theta} \\
            0
        \end{pmatrix}.        
    \end{equation}
    Now, differentiating~\eqref{eq:v-n-change} any fixed number $r$ of times with respect to $\alpha$, we get in the right hand side $O(n^r)$ terms (each of $r$ derivatives can land on any of $n+1$ factors). At the same time, differentiating $r'$ times in $\alpha$ the factor $\kappa^s_{\alpha,\theta + j \alpha}$, $0\le j <n$, one gets a factor that is $O(n^{r'})$, hence each of the such terms does not exceed $C n^r \rho^n$, and thus in total $r$-th derivative does not exceed $C' n^{2r} \rho^n$, that is exponentially decreasing in $n$ for any fixed~$r$.

    The exponential decay of the derivatives for $n<0$ is handled in the same way. 
    %Finally, if we could not choose vectors in the stable/unstable direction in the continuous way, it means that the corresponding monodromy of the stable direction leads to the direction-reversal; we then pass to the twofold cover of the base (in the same way as it is done in \cite{HA}), and repeat the above arguments.

\end{proof}

\begin{proof}[Proof of Proposition~\ref{p:Phi-bounds}]
From assumption, the function $\Phi$ can be extended as a holomorphic function to some complex neighbourhood $U_{\C}$ of~$\Spl_{\alpha,\lambda}$. Choosing a smaller neighbourhood $U_1$ if necessary, we can assume (see Fig.~\ref{fig:contour}) that~$U_1$ has a smooth boundary~$\gamma:=\partial U_1$.

\begin{figure}
    \centering
    \includegraphics[width=0.7\linewidth]{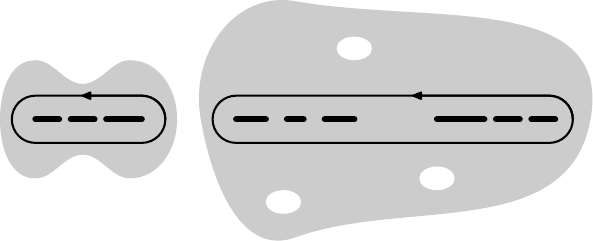}
    \caption{The spectrum $\Sigma^+_{\alpha_0,\lambda_0}$, domain $U_{\C}$ and the boundary (and the integration contour)~$\partial U_1$.}
    \label{fig:contour}
\end{figure}

Now, applying the Cauchy's formula to the function $\Phi$, we can represent 
\[
    \Phi(E) = \frac{1}{2\pi i} \int_{\gamma} \frac{\Phi(z) \, dz}{z-E}.
\]
Then, using the Cauchy formula for operator valued  functions (e.g. see \cite[Appendix A]{PT}), one can write  %for every $(\alpha,\lambda)\in B$ 
\begin{equation}\label{eq:Green}
    \Phi(H_{\alpha,\lambda,\theta}) = \frac{1}{2\pi i} \int_{\partial U_1} \Phi(z) \cdot (-G_{\alpha,\lambda,\theta}(z)) \, dz,   
\end{equation}
where $G_{\alpha,\lambda,\theta}(z)$ is the Green's function, defined by~\eqref{eq:Green-def}. 

For every point $z\in\gamma$ of the integration contour, there exists a neighbourhood $U_z^{\alpha}\times U_z^{\lambda}\times U_z^{\gamma}$ of $(\alpha,\lambda,z)$, at which the uniform estimates of Lemma~\ref{l:Green} hold with some $C_{1,z},c_{1,z}$. Extracting a finite cover of $\gamma$ by some $U^{\gamma}_{z_r}$ and integrating these estimates elementwise, we obtain the desired decay bounds with 
\[
c_0= \min_r c_{1,z_r}, \quad C_0=\frac{|\gamma| \cdot \max_{\gamma} |\Phi(z)|}{2\pi} \cdot \max_r C_{1,z_r}, 
\]
that hold in the neighbourhood $\bigcap_j (U_{z_j}^{\alpha} \times U_{z_j}^{\lambda})$ of~$(\alpha,\lambda)$.

 Finally, the analyticity with respect to $\lambda$ and smoothness with respect to $\alpha$ follow from the same properties of the elements of the Green function and the uniform estimates from Lemma~\ref{l:Green}.
\end{proof}

%\vspace{1cm}

\subsection{Four traces formula}\label{s:four-traces}
 %: Proposition~\ref{p.A-A-duality}

The estimates of Proposition~\ref{p:Phi-bounds} state that in its assumptions, the elements of $(\Phi(H_{\alpha,\lambda,\theta}))_{jl}$ depend analytically on $\lambda$ and decay exponentially in~$|j-l|$. Hence, we can define the traces as the sums of anti-diagonal elements:
\begin{equation}\label{eq:tr-Phi-0}
\tr \left(\Phi(H_{\alpha,\lambda,0}) R_0\right) := \sum_{j\in \Z} (\Phi(H_{\alpha,\lambda,0}))_{j,-j},    
\end{equation}
\begin{equation}\label{eq:tr-Phi-fa}
\tr \left(\Phi(H_{\alpha,\lambda,\fa}) R_{\fa} \right) := \sum_{j\in \Z} (\Phi(H_{\alpha,\lambda,\fa}))_{j,-1-j},    
\end{equation}
\begin{equation}\label{eq:tr-Phi-1-2}
\tr \left(\Phi(H_{\alpha,\lambda,\frac{1}{2}}) R_{\frac{1}{2}} \right) := \sum_{j\in \Z} (\Phi(H_{\alpha,\lambda,\frac{1}{2}})))_{j,-j},    
\end{equation}
\begin{equation}\label{eq:tr-Phi-1-a-2}
\tr \left(\Phi(H_{\alpha,\lambda,\frac{1+\alpha}{2}}) R_{\frac{1+\alpha}{2}} \right) := 
\sum_{j\in \Z} (\Phi(H_{\alpha,\lambda,\frac{1+\alpha}{2}}))_{j,-1-j},    
\end{equation}

\begin{theorem}[Four-traces formula]\label{t:four-traces}
% \todo[inline]{The sum of the four traces}
For \emph{any} $\alpha$, any $\lambda<1$, and any function $\Phi$ that is analytic in some neighbourhood of $\Spl_{\alpha,\lambda}$, one has
    \begin{multline}\label{eq:four-traces}
        \int_{\Si_{\alpha,\lambda}} \Phi'(E) \, d\mmu_{\alpha,\lambda}(E) = 
        \tr \Phi(H_{\alpha,\lambda,\frac{1+\alpha}{2}})R_{\frac{1+\alpha}{2}} + \tr \Phi(H_{\alpha,\lambda,\fa}) R_{\fa} +
        \\
        +\tr \Phi(H_{\alpha,\lambda,\frac{1}{2}}) R_{\frac{1}{2}}- \tr \Phi(H_{\alpha,\lambda,0}) R_0,
    \end{multline}
    where the traces in the right hand side are defined by~\eqref{eq:tr-Phi-0}--\eqref{eq:tr-Phi-1-a-2}.
    For $\lambda>1$, the integral differs from the right hand side of~\eqref{eq:four-traces} by a change of the sign.
\end{theorem}

\begin{proof}
    First, assume that $\alpha=\frac{p}{q}$ is rational. We will need the following lemma.
    \begin{lemma}\label{l:projection}
        For any $\alpha=\frac{p}{q}$, $\lambda$, $\theta$ and any function $\Phi$ that is analytic in some neighbourhood of $\Spl_{\alpha,\lambda}$, the elements of $\Phi(\Hfp_{\alpha,\lambda,\theta})$ are given by
        \begin{equation}\label{eq:Phi-per-proj}
            (\Phi(\Hfp_{q,\theta}))_{j,l} = \sum_{m\in\Z} (\Phi(H_{\alpha,\lambda,\theta}))_{j+mq,l}
        \end{equation}
        In the same way, the elements of $\Phi(\Hfa_{\alpha,\lambda,\theta})$
        are given by
        \begin{equation}\label{eq:Phi-anti-proj}
            (\Phi(\Hfa_{q,\theta}))_{j,l} = \sum_{m\in\Z} (-1)^m (\Phi(H_{\alpha,\lambda,\theta}))_{j+mq,l}.
        \end{equation}
    \end{lemma}
         Informally (and completely non-rigorously!) speaking, the statement of the lemma could be thought as a statement that the ``projections'' $\sum_{m\in \Z} \mT_q^m$ and $\sum_{m\in \Z} (-1)^m \mT_q^m$ to the $\mT_q$-generalized eigenspaces $\Vp_q$ and $\Va_q$ respectively commute with~$H_{\alpha,\lambda,\theta}$ and thus with $\Phi(H_{\alpha,\lambda,\theta})$. However, these operator series diverge, so for a formal proof we will need to find a workaround.

    \begin{proof}[Proof of Lemma~\ref{l:projection}]
         Consider first a particular case $\Phi(E)=\frac{1}{E-z}$, where $z\notin \Spl_{\alpha,\lambda}$. Note that due to Lemma~\ref{l:Green}, for every $j,l$ the series in the right hand side of~\eqref{eq:Phi-per-proj} converges exponentially. Hence, for a fixed $l$ this series defines a vector $v_l\in \Vp_q$,
        \begin{equation}\label{eq:v-j-def}
            (v_l)_j = \sum_{m\in\Z} (\mT_q^m (H_{\alpha,\lambda,\theta}-z)^{-1} \ve_l)_{j}, %\quad 
        \end{equation}
        where $(\ve_l)_i = \delta_{li}$ are the base vectors of~$\ell^2(\Z)$.
        Now, apply to $v_l\in\Vp_q$ the operator $\Hfp_{q,\theta}-z$. As any $j$-th coordinate of its image depends only on the coordinates~$(v_l)_{j-1},(v_l)_{j},(v_l)_{j+1}$, each of which is given by a convergent series~\eqref{eq:v-j-def}, we have 
        \begin{multline}            
            ((\Hfp_{q,\theta}-z)v_l)_j = \sum_m \left( (\Hfp_{q,\theta}-z) (\mT_q^m (H_{\alpha,\lambda,\theta}-z)^{-1} \ve_l) \right)_j
            \\
            = \sum_m \left( \mT_q^m \ve_j \right)_l = \begin{cases}
                1 & j \equiv l \mod q \\
                0 & \text{otherwise}.
            \end{cases}
        \end{multline}
        Hence, $(\Hfp_{q,\theta}-z)v_l=e_l\in \Vp_q$, and thus $v_l=(\Hfp_{q,\theta}-z)^{-1} e_l= \Phi(\Hfp_{q,\theta}) e_l$, establishing~\eqref{eq:Phi-per-proj} for this particular function~$\Phi$. Now, in the same way as in the proof of Proposition~\ref{p:Phi-bounds}, we can apply Cauchy's formula
        \[
            \Phi(E) = \frac{1}{2\pi i} \int_{\gamma} \frac{\Phi(z) \, dz}{z-E},
        \]
        and represent the operators $\Phi(\Hfp_{q,\theta})$ and $\Phi(H_{\alpha,\lambda,\theta})$ in both sides of~\eqref{eq:Phi-per-proj} as contour integrals. The uniform exponential decay estimates of Proposition~\ref{p:Phi-bounds} allow to interchange (for any fixed $j,l$) the integration over $z$ and summation over $m$, thus implying
        \begin{multline}
                (\Phi(\Hfp_{q,\theta}))_{j,l} = \frac{1}{2\pi i} \int_{\gamma} \Phi(z) ((\Hfp_{q,\theta}-z)^{-1})_{j,l} \, dz 
                \\
                =\frac{1}{2\pi i} \int_{\gamma} \Phi(z) \bigl(\sum_{m\in \Z} (H_{\alpha,\lambda,\theta}-z)^{-1}\bigr)_{j+mq,l} \, dz 
                \\
                = \sum_{m\in \Z}  \frac{1}{2\pi i} \int_{\gamma} \Phi(z) \bigl((H_{\alpha,\lambda,\theta}-z)^{-1} \bigr)_{j+mq,l} \, dz 
                = \sum_{m\in \Z}  (\Phi(H_{\alpha,\lambda,\theta}))_{j+mq,l},
        \end{multline}
        what concludes the proof for~\eqref{eq:Phi-per-proj}. The equality~\eqref{eq:Phi-anti-proj} is proven in exactly the same way, considering the sums with an additional $(-1)^m$ factor.
    \end{proof}
        %  {Let us first present and informal (and non-rigorous!) argument. Note that the left hand side of~\eqref{eq:Phi-per-proj} is the $l$-th coorinate of  $\Phi(H_{\alpha,\lambda,\theta})$-image of the base vector $e_j\in \Vp$, defined by~\eqref{eq:e-base},
        % and this vector can \emph{informally} be considered as a sum of translates of the coordinate one, $\mathbf{e}_j$: 
        % \[
        %     e_j = \sum_{m\in\Z} \mT_q^m, \quad (\mathbf{e}_j)_n = \delta_{jn}.
        % \]
        % Meanwhile, the right hand side is (again, informally speaking) the $l$-th coordinate of the image $\sum_m \mT_q^{m} \left( \Phi(H_{\alpha,\lambda,\theta}) e_j \right)$. Now, as $\mT_q$ and $H_{\alpha,\lambda,\theta}$ commute (recall that $\alpha=\frac{p}{q}$), so would the (not well-defined) sum $\sum_{m\in \Z} \mT_q^m$ and $\Phi(H_{\alpha,\lambda,\theta})$; applying them to~$\mathbf{e}_j$ in different order would imply the desired~\eqref{eq:Phi-per-proj}. Analogously, replacing $\sum_m \mT_q^m$ by $\sum_m (-1)^m \mT_q^m$, one would get~\eqref{eq:Phi-anti-proj}. Unfortunately, the series appearing in these arguments do not converge in $\ell^2$ or $\ell^{\infty}$, and thus have no formal meaning; the proof thus requires additional arguments.}       
%        \todo[inline]{Projections on the generalized eigenspaces for the translation operator $\mT_q$ that commutes with $H_{\alpha,\lambda,\theta}$}

    Lemma~\ref{l:projection} allows to rewrite the right hand side of~\eqref{eq:Phi-trace} and of~\eqref{eq:int-phi-super}:
    \begin{lemma}\label{l:2-to-4}
        For any $\alpha$, $\lambda\in \R$, and any function $\Phi$ that is analytic in some neighbourhood of $\Spl_{\alpha,\lambda}$, one has
        \begin{multline}\label{eq:2-to-4}
            \tr \left(\Phi(\Hfp_{q,\ta}) R_{\ta}\right) - \tr \left(\Phi(\Hfa_{q,0}) R_0\right) = 
            \tr \left( \Phi(H_{\alpha,\lambda,\fa}) R_{\fa}\right) +
            \\
            + \tr \left(\Phi(H_{\alpha,\lambda,\frac{1+\alpha}{2}})R_{\frac{1+\alpha}{2}}\right) +\tr \left(\Phi(H_{\alpha,\lambda,\frac{1}{2}}) R_{\frac{1}{2}} \right) - \tr \left( \Phi(H_{\alpha,\lambda,0}) R_0 \right).
        \end{multline}
    \end{lemma}
    \begin{proof}
    The first of the traces in the left hand side of~\eqref{eq:2-to-4} can be written using Lemma~\ref{l:projection} as 
    \[
            \tr \left(\Phi(\Hfp_{q,\ta}) R_{\ta}\right) = \sum_{j=1}^{q} \left(\Phi(\Hfp_{q,\ta}) \right)_{j, -1-j} = \sum_{j=1}^q \sum_{m\in\Z} \Phi(H_{\alpha,\lambda,\ta})_{j+mq,-1-j}.
    \]
    Now, splitting the summation into two sums over even and odd values of~$m$, one gets 
    \begin{multline}\label{eq:transforming}        
        \sum_{j=1}^q \sum_{m\in\Z} \Phi(H_{\alpha,\lambda,\ta})_{j+mq,-1-j}=
        \\
        =\sum_{j=1}^q \sum_{m\in\Z} \Phi(H_{\alpha,\lambda,\ta})_{j+2mq,-1-j} + \sum_{j=1}^q \sum_m \Phi(H_{\alpha,\lambda,\ta})_{j+(2m-1)q,-1-j}
        \\
        = \sum_{j=1}^q \sum_{m\in\Z} \Phi(H_{\alpha,\lambda,\ta})_{j+mq,-1-(j+mq)} + \sum_{j=1}^q \sum_m \Phi(H_{\alpha,\lambda,\ta})_{j+mq ,-1+q-(j+mq)}
        \\ 
        = \sum_{j'\in\Z} \Phi(H_{\alpha,\lambda,\ta})_{j',-1-j'} +
        \sum_{j'\in\Z} \Phi(H_{\alpha,\lambda,\ta})_{j',-1+q-j'},
    \end{multline}
    as the new index $j'=j+mq$ takes every integer value exactly once. Now, as $\ta\in \left\{\frac{\alpha}{2}, \frac{1+\alpha}{2}\right\}$, the first of two sums in the right hand side of~\eqref{eq:transforming} by definition~(\eqref{eq:tr-Phi-fa} or~\eqref{eq:tr-Phi-1-a-2}) is equal to $\tr \left(\Phi(H_{\alpha,\lambda,\ta}) R_{\ta}\right)$. For the second sum in the right hand side of~\eqref{eq:transforming}, note that for every phase $\theta$ and every~$i,j,r$,
    \[
        \Phi(H_{\alpha,\lambda,\theta})_{i+r,j+r}=\Phi(H_{\alpha,\lambda,\theta+r\alpha})_{i,j}.
    \]
    Take $r:=\lfloor q/2\rfloor$ and let $j':=j''+r$, then 
    \begin{multline}\label{eq:trace-theta-prime}
        \sum_{j'\in\Z} \Phi(H_{\alpha,\lambda,\ta})_{j',-1+q-j'} = \sum_{j''\in\Z} \Phi(H_{\alpha,\lambda,\ta})_{j''+r,-1+q-(j''+r)}
        \\
        = \sum_{j''\in\Z} \Phi(H_{\alpha,\lambda,\ta+r\alpha})_{j'',(-1+q-2r)-j''}.
    \end{multline}
    Now, the phase $\theta':=\ta+r\alpha$ is the second point where the orbit of $\ta$ under the rotation by $\alpha=\frac{p}{q}$ intersects the set $\Th_{\alpha}$: for even $q$, one has 
    \[
    \theta'=\ta + \frac{q}{2}\cdot \frac{p}{q} \equiv \ta + \frac{1}{2} \mod 1,
    \]
    while for odd $q$,
    \[
    \theta' = \ta + \frac{q-1}{2} \cdot \frac{p}{q} \equiv \frac{1}{2} \mod 1.
    \]
    The value $-1+q-2r$ is equal to $-1$ and $0$ for even and odd $q$ respectively, and in both cases the right hand side of~\eqref{eq:trace-theta-prime} is thus by definition equal to
    \[
        \tr \left(\Phi(H_{\alpha,\lambda,\theta'}) R_{\theta'}\right).
    \]
    Adding up, we see that the first of the traces in the left hand side of~\eqref{eq:2-to-4} is equal to 
    \begin{equation}\label{eq:tr-ta-final}
        \tr \left(\Phi(\Hfp_{q,\ta}) R_{\ta}\right) = \tr \left(\Phi(H_{\alpha,\lambda,\ta}) R_{\ta}\right)+ \tr \left(\Phi(H_{\alpha,\lambda,\theta'}) R_{\theta'}\right).        
    \end{equation}

    In the same way, the second trace in the left hand side of~\eqref{eq:2-to-4} equals
     \begin{multline}\label{eq:transforming-2}        
        \tr \left(\Phi(\Hfa_{q,0}) R_{0}\right) = \sum_{j=1}^{q} \left(\Phi(\Hfp_{q,0}) \right)_{j, -j} = \sum_{j=1}^q \sum_m (-1)^m \Phi(H_{\alpha,\lambda,0})_{j+mq,-j}
        \\
        =\sum_{j=1}^q \sum_m \Phi(H_{\alpha,\lambda,0})_{j+2mq,-j} - \sum_{j=1}^q \sum_m \Phi(H_{\alpha,\lambda,0})_{j+(2m-1)q,-j}
        \\
        = \sum_{j=1}^q \sum_m \Phi(H_{\alpha,\lambda,\ta})_{j+mq,-(j+mq)} - \sum_{j=1}^q \sum_m \Phi(H_{\alpha,\lambda,0})_{j+mq ,q-(j+mq)}
        \\ 
        = \sum_{j'} \Phi(H_{\alpha,\lambda,0})_{j',-j'} -
        \sum_{j'} \Phi(H_{\alpha,\lambda,0})_{j',q-j'}
        \\ 
        = \tr \left(\Phi(H_{\alpha,\lambda,0}) R_{0}\right) - 
        \tr \left(\Phi(H_{\alpha,\lambda,\theta''}) R_{\theta''}\right),
    \end{multline}   
    where we again use an index shift by $r':=\lceil q/2 \rceil$, and $\theta''=r'\alpha$ is the fourth of the angles of~$\Th_{\alpha}$ (the second one on the orbit of~$0$). 
    
    Subtracting~\eqref{eq:transforming-2} from~\eqref{eq:tr-ta-final} completes the proof of Lemma~\ref{l:2-to-4}.
    \end{proof}
    Now we are ready to complete the proof of Theorem~\ref{t:four-traces}.
    Joining Lemma~\ref{l:2-to-4} with Propositions~\ref{p:Phi-tr} and~\ref{p:Phi-tr-super}, we immediately get the conclusion of Theorem~\ref{t:four-traces} for rational $\alpha=\frac{p}{q}$ (for $\lambda<1$ and $\lambda>1$ respectively).

    Next, let $\ba$ be irrational and assume that for some $\bl$ the function $\Phi$ is analytic in a neighbourhood of $\Spl_{\ba,\bl}$. Denote by $\mI_{\Phi}(\alpha,\lambda)$ the right hand side of~\eqref{eq:2-to-4}; Proposition~\ref{p:Phi-bounds} implies that $\mI_{\Phi}(\alpha,\lambda)$ is continuous in some neigbourhood of $(\ba,\bl)$.

    If $\bl<1$, a result from~\cite{JM12} that we have already mentioned in the introduction, states that there exists a sequence $\alpha_n=\frac{p_n}{q_n}$ that converges to $\ba$, such that $\Leb(\Si_{\alpha_n,\bl}\SymDiff \Si_{\ba,\bl})\to 0$. We thus have 
    \begin{multline}\label{eq:four-proof-sub}
        \int \Phi'(E) \, d\mmu_{\ba,\bl}(E) = \lim_{n\to\infty }\int \Phi'(E) \, d\mmu_{\alpha_n,\bl}(E) =
        \\
        =\lim_{n\to\infty } \mI_{\Phi}(\alpha_n,\bl) = \mI_{\Phi}(\ba,\bl).
    \end{multline}
    Here the second equality holds as we have already established conclusion~\eqref{eq:four-traces} for rational values $\alpha_n$ of the frequency, and the last one follows from the continuity of the right hand side $\mI_{\Phi}(\alpha,\lambda)$. We have thus proven~\eqref{eq:four-traces} for $\lambda=\bl<1$.

    On the other hand, if $\bl>1$, then we can take the sequence $\alpha_n\to \ba$ such that  
    $\Leb(\Si_{\alpha_n,\frac{1}{\bl}}\Delta \Si_{\ba,\frac{1}{\bl}})\to 0$. Then we also have
    \[
        \Leb(\Si_{\alpha_n,\bl}\SymDiff \Si_{\ba,\bl}) =
        \Leb\left(\bl\Si_{\alpha_n,\frac{1}{\bl}}\SymDiff \bl\Si_{\ba,\frac{1}{\bl}}\right) \to 0,
    \]
    and repeating~\eqref{eq:four-proof-sub} (with a sign change coming from $\bl>1$) completes the proof of Theorem~\ref{t:four-traces} for $\lambda=\bl>1$.
\end{proof}

\section{Proofs of main results}\label{s:proofs}

\begin{proof}[Proof of Theorem~\ref{t.anmu}]
    Theorem~\ref{t.anmu} immediately follows from Theorem~\ref{t:four-traces}, providing the trace formula for this integral, and Proposition~\ref{p:Phi-bounds}, implying the smoothness in $\alpha$ and analyticity in $\lambda$ of these traces.
\end{proof}

\begin{proof}[Proof of Theorem~\ref{t:poly}]    
    Theorem~\ref{t:poly} also follows from Theorem~\ref{t:four-traces}, applied for $\Phi(E)=E^{2k+1}$. It suffices to check that each of the four traces appearing in the right hand side of~\eqref{eq:four-traces} has the required form. 

    First, note that the four-traces formula~\eqref{eq:four-traces}, applied to the integral of $\Phi'(E)=(2k+1) E^{2k}$, can be simplified due to the symmetry $\theta\mapsto \frac{1}{2}+\theta$, that changes the signs of all the diagonal elements. The sum in the right hand side then has the form 

% \begin{coro}\label{c:simple}
%     The right hand side of~\eqref{eq:moment-traces} equals 
    \begin{multline}\label{eq:simple}
        2 \left( \tr H_{\alpha,\lambda,\ta}^{2k+1} R_{\ta} - \tr H_{\alpha,\lambda,0}^{2k+1} R_0 \right)  \\
        = 2 \left( \sum_{i=-k-1}^{k} \left(H_{\alpha,\lambda,\ta}^{2k+1}\right)_{i,-1-i} - \sum_{i=-k-1}^{k+1} \left(H_{\alpha,\lambda,0}^{2k+1}\right)_{i,-i} \right)
    \end{multline}    
% \end{coro}

% Now, take $P_{2k}(\alpha,\lambda)$ to be $\frac{1}{2k+1}$ times the right hand side of~\eqref{eq:simple}; the coefficient is due to the fact  that the anti-derivative of $\varphi(E)=E^{2k}$ is $\Phi(E)=\frac{1}{2k+1}E^{2k+1}$. Then, $P_{2k}(\alpha,\lambda)$ is indeed of the form~\eqref{eq:P}, a polynomial in $\lambda$ of degree $2k+1$ with the coefficients that are trigonometric polynomials in $\alpha$ and leading term $(-\frac{2^{2k+2}}{2k+1})$. 

Now, the right hand side of~\eqref{eq:simple} is a sum of a finite number of matrix elements that are polynomial in $\lambda$ (see Eqs.~\eqref{eq:H-R-0}--\eqref{eq:H-R-12}). The form of the leading term $-2^{2k+2}$ of the polynomial is also immediate: it comes from the term $-V(0)^{2k+1}=-2^{2k+1}\lambda^{2k+1}$ in the second sum. The fact that the coefficients of this polynomial are trigonometric polynomials in $\alpha$ is easy to see for the second sum. For the first one, the diagonal elements of $H_{\alpha,\lambda,\ta}$ are all of the form $2\lambda \cdot \cos (2j+1)\pi\alpha$, and the cosine is given by an odd (Chebyshev) polynomial of $\cos \pi \alpha$. Now, it is not difficult to see that all the contributions to the first summand come from the product of an even number of such terms. Such a product is thus an even power of~$\lambda$, multiplied by an even polynomial of $\cos \pi \alpha$. As $\cos^2 \pi\alpha =\frac{1+\cos 2\pi\alpha}{2}$, the latter is also a trigonometric polynomial of $\cos 2\pi\alpha$.
\end{proof}

\begin{proof}[Proof of Theorem~\ref{t.Leb-cont}]    
Theorem~\ref{t.Leb-cont} follows from Theorem~\ref{t:poly}, as for measures of bounded mass, supported on a compact interval, convergence of all moments implies weak convergence, e.g. see~\cite[Theorem~3.3.26]{Du}. We reproduce here the corres\-ponding argument for completeness.

Specifically, due to Theorem~\ref{t:poly}, the moments $c_k(\alpha,\mu)$ depend continuously on~$\alpha$ and~$\lambda$. Hence, the same applies to the dependence on $(\alpha,\lambda)$ of integrals 
\[
    \int g(E) \, d\mmu_{\alpha,\lambda}
\]
of any polynomial~$g(E)$.

For $\lambda<1$, all the measures $\mmu_{\alpha,\lambda}$ are supported inside the interval $[-4,4]$ and do not exceed $8$ in total mass. Now, any continuous function $f$ on this (compact) interval can be uniformly approximated by polynomials, thus implying that its integrals w.r.t. $\mmu_{\alpha,\lambda}$ depend continuously on $(\alpha,\lambda)$ as well (as uniform limits of continuous functions). This by definition implies continuous dependence of $\mmu_{\alpha,\lambda}$ on~$(\alpha,\lambda)$.

Finally, to handle the case $\lambda>1$ it suffices to note that $\mmu_{\alpha,\frac{1}{\lambda}}$ is transformed to~$\mmu_{\alpha,\lambda}$ by scaling and multiplication by $\lambda$, and both of these operations are continuous.
\end{proof}
%\todo[inline]{Or should we say that for every $\lambda\le C$ we have compact interval and upper bound for the mass?}
% \begin{proof}[Proof of ]    
% \end{proof}
\begin{proof}[Proof of Theorem~\ref{t.analytic}]    
As we have already seen in Remark~\ref{r.particular}, Theorem~\ref{t.analytic} is a particular case of Theorem~\ref{t.anmu} for the function $\varphi$ that is equal to~$1$ in a neighbourhood $U_1$ of the region between the two gaps, and $0$ in a neighbourhood $U_2$ of the rest of the spectrum.  
\end{proof}

\begin{proof}[Proof of Theorem~\ref{t.anmutilde}]    
For any $\alpha$ and $\lambda\neq 1$, define
\begin{equation}\label{eq:mu-def}
    \tmu_{\alpha,\lambda}:=\frac{1}{|4-4\lambda|} \mmu_{\alpha,\lambda}.
\end{equation}
For a function $\Phi$ that is analytic in a neighbourhood of $\Spl_{\alpha,\lambda}$, recall that $\mI_{\Phi}(\alpha,\lambda)$ is defined as the right hand side of~\eqref{eq:four-traces}, and this function is then smooth in $\alpha$ and analytic in $\lambda$ in the neighbourhood of~$(\alpha_0,\lambda_0)$. Then, from Theorem~\ref{t:four-traces},
\begin{multline}\label{eq:Phi-tmu}
    \int \Phi'(E) \, d\tmu_{\alpha,\lambda}(E) = \frac{1}{4\cdot |1-\lambda|} \int \Phi'(E) \, d\tmu_{\alpha,\lambda}(E) 
    \\ = \begin{cases}
        \frac{1}{4(1-\lambda)} \mI_{\Phi}(\alpha,\lambda), & 0<\lambda <1\\
        -\frac{1}{4\cdot |1-\lambda|} \mI_{\Phi}(\alpha,\lambda), & \lambda >1
    \end{cases}
    \\ = \frac{1}{4(1-\lambda)} \mI_{\Phi}(\alpha,\lambda).
\end{multline}
This implies that the integral coincides with a function that is smooth in $\alpha$ and analytic in~$\lambda$ away from~$\lambda=1$.

Now, assume $\Phi$ is analytic in in a neighbourhood of $\Spl_{\alpha,1}$. Then, the measure $\Leb(\Si_{\alpha,\lambda})$ tend to zero as $\lambda\to 1$, hence so does the integral $\int_{\Si_{\alpha,\lambda}} \Phi(E) d\mmu_{\alpha,\lambda}(E)$ and thus the function $\mI_{\Phi}(\alpha,\lambda)$. At the same time, in this assumption this function is analytic in $\lambda$ in a neighbourhood of $(\alpha,1)$, hence $\mI_{\Phi}(\alpha,1)=0$ and thus the quotient 
\[
\widetilde{\mI}_{\Phi} (\alpha,\lambda) := \frac{1}{4(1-\lambda)} \mI_{\Phi}(\alpha,\lambda)
\]
can be extended continuously to $\lambda=1$, and the result of this continuation is also analytic in $\lambda$ and smooth in $\alpha$.
%\todo[inline]{Should we recall here the trick $\frac{g(x)}{x} = \int_0^1 g'(tx) dt$? Reference to Hadamard's Lemma?}

Finally, define $\tmu_{\alpha,1}$ as any accumulation point of the family of measures $\tmu_{\alpha,1-\frac{1}{n}}$. Then all the moments of the measures $\tmu_{\alpha,\lambda}$ converge to the moments of~$\tmu_{\alpha,1}$ due to the continuity of $\widetilde{\mI}_{E^{k+1}}(\alpha,\lambda)$. Thus, the same applies to the integral of any continuous function (that can be uniformly approximated by polynomials on a large segment of the real line). Thus, the family of measures $\tmu_{\alpha,\lambda}$ is continuous at $\lambda=1$, and for every function $\Phi$ that is analytic in the neighbourhood of $\Spl_{\alpha,1}$, one has 
\[
\int \Phi'(E) \, d\tmu_{\alpha,\lambda} = \widetilde{I}_{\Phi}(\alpha,\lambda)
\]
for any $\lambda>0$, including $\lambda=1$.
\end{proof}

\begin{proof}[Proof of Proposition~\ref{p.A-A-duality}]
We have 
\[
P_{2k}(\alpha,\lambda)= \frac{2}{2k+1} (\tr H_{\alpha,\lambda,\frac{\alpha}{2}}^{2k+1} R_{\frac{\alpha}{2}} - \tr H_{\alpha,\lambda,0}^{2k+1} R_{0}),
\]
both summands in the right hand side being polynomial in $\lambda$. As Eq.~\eqref{eq:A-A-duality-bis} in Example~\ref{ex:compute-second} states that
\[
\tr H_{\alpha,\lambda,\frac{\alpha}{2}}^{2k+1} R_{\frac{\alpha}{2}} = - \lambda^{2k+1} \tr H_{\alpha,\frac{1}{\lambda},0}^{2k+1} R_{0},
\]
the conclusion of the proposition follows.

% We know from Theorems~\ref{t.anmutilde} and~\ref{t:four-traces} that the moment
% \[
%     \tc_{2k}(\alpha,\lambda)=\int E^{2k} \, d\tmu_{\alpha,\lambda}
% \]
% is a polynomial function of~$\lambda$. At the same time, the measure $\tmu_{\alpha,\lambda}$ is a push-forward of  $\tmu_{\alpha,\frac{1}{\lambda}}$ by rescaling by~$\lambda$, hence 
% \[
%     \tc_{2k}(\alpha,\lambda) = \lambda^{2k} \, \tc_{2k}(\alpha,\frac{1}{\lambda}).
% \]
% Returning to the moments of non-normalized measures $\mmu_{\alpha,\lambda}=|4-4\lambda| \tmu_{\alpha,\lambda}$, we get 
% \[
% c_{2k}(\alpha,\lambda) = (4-4\lambda) \, \tc_{2k}(\alpha,\lambda), \quad \lambda<1,
% \]
% and thus for $\lambda>1$
% \begin{multline}
% c_{2k}(\alpha,\lambda) = (4\lambda-4) \tc_{2k}(\alpha,\lambda) = (4\lambda-4) \lambda^{2k} \, \tc_{2k}(\alpha,\frac{1}{\lambda}) 
% \\
% =  (4\lambda-4)\lambda^{2k} \cdot \frac{c_{2k}(\alpha,\frac{1}{\lambda})}{4-4\frac{1}{\lambda}} 
% = \lambda^{2k+1} c_{2k}(\alpha,\frac{1}{\lambda}).
% \end{multline}

% At the same time, we know from Theorem~\ref{t:four-traces} that polynomial function $P_{2k}(\alpha,\lambda)$, given by the right hand side of~\eqref{eq:four-traces} for $\lambda<1$, differs from $c_{2k}(\alpha,\lambda)$ by the sign for $\lambda>1$.  Hence, for $\lambda>1$
% \[
% P_{2k}(\alpha,\lambda) = -c_{2k} (\alpha,\lambda) = -\lambda^{2k+1} c_{2k}(\alpha,\frac{1}{\lambda}) = - \lambda^{2k+1} P_{2k}(\alpha,\frac{1}{\lambda}).
% \]
% This completes the proof of Proposition~\ref{p.A-A-duality}.
\end{proof}

\section{Non-analytic dependence on the angle~$\alpha$}\label{s:ex-smooth}

In the conclusions of Theorems~\ref{t.analytic},  \ref{t.anmu},~\ref{t.anmutilde} we only state $C^{\infty}$-smoothness in~$\alpha$. We believe that, contrary to the polynomial case of Theorem~\ref{t:poly}, this dependence is actually \emph{not} analytic, and in this section, we present arguments in this direction. 

Namely, the proof of Theorem~\ref{t.anmu} uses the decomposition $\C^2=E^s_x \oplus E^u_x$ into the stable and unstable directions; its existence (and analyticity in $x,E$ and $\lambda$) is proven via the graph transform principle. The example below shows that a graph constructed in this way may depend on~$\alpha$ in a non-analytic way even in a simpler setting.

Namely, consider the following one-dimensional skew product over the rotation by $\alpha$: 
\[
F_{\alpha}(x,y) = (x+\alpha, \frac{1}{2} (y + \varphi(x))), \quad (x,y)\in \Sc\times \R,
\]
where $\varphi(x)$ is a given analytic function on the circle. The same graph transform principle implies that this skew product admits an invariant section $y=\psi_{\alpha}(x)$. Moreover, writing down the formula for the $n$-th image of zero section, and passing to the limit, one gets an explicit formula for this section: 
\begin{equation}\label{eq:psi}
\psi_{\alpha}(x) = \sum_{n=1}^{\infty} \frac{1}{2^n} \, \varphi(x-n\alpha).
\end{equation}

Now, on the one hand, the function $\psi_{\alpha}(x)$ is infinitely differentiable in $\alpha$, and its $k$-th derivative is given by termwise differentiation of the series~\eqref{eq:psi}:
\begin{equation}\label{eq:k-th-derivative}
\frac{\partial^k}{\partial \alpha^k}\psi_{\alpha}(x) = (-1)^k \sum_{n=1}^{\infty} \frac{n^k}{2^n} \, \cdot \varphi^{(k)}(x-n\alpha).    
\end{equation}
On the other, for large $k$ the elements of the series~\eqref{eq:k-th-derivative} seem to be too large to add up with the analyticity in~$\alpha$. Namely, dividing \eqref{eq:k-th-derivative} by $k!$ and considering the absolute value of the term with $n=k$, one gets
\[
\frac{k^k}{2^k} \cdot \frac{|\varphi^{(k)} (x-k\alpha)|}{k!}.
\]
For an analytic function $\varphi$, the second factor is allowed to grow exponentially, and then this particular term has a superexponential growth in~$k$, hinting a possible non-analyticity. Though the above arguments are not a formal proof (large terms might cancel out), the proposition below shows that generically the dependence on~$\alpha$ is indeed non-analytic. To state it, let us introduce some notations. Namely, let~$\Psi$ be the linear map that associates to a function $\varphi(x)$ the right hand side $\psi_{\alpha}(x)$ of~\eqref{eq:psi}, and let 
\[
    \Ta_k(\varphi,x,\alpha_0) = \frac{1}{k!}  \left.
    \frac{\partial^k}{\partial \alpha^k}\right|_{\alpha=\alpha_0}\Psi(\varphi)(x,\alpha)
\]
be the Taylor coefficients of $\psi_{\alpha}$ in the variable~$\alpha$ at the point $(x,\alpha_0)$.

\begin{prop}\label{p:non-analytic}
There exists an analytic function $\varphi$ on the circle, for which for the Taylor series in $\alpha$ of the corresponding function~$\psi_{\alpha}=\Psi(\varphi)$, defined by~\eqref{eq:psi}, has zero convergence radius at every point $(x,\alpha_0)$. Namely,
\begin{equation}\label{eq:c-growth}
\forall x\in \Sc, \quad \forall \alpha_0\in \Sc \quad \limsup_k \sqrt[k]{|\Ta_k(\varphi,x,\alpha_0)|} = +\infty.
\end{equation}
In particular, the dependence of $\psi_{\alpha}(x)$ on $\alpha$ is not analytic at any point $(x,\alpha_0)$.
\end{prop}
\begin{coro}
    For any $(x,\alpha_0)$ and any function $\phi_0(x)$ on the circle, consider a family of functions $\phi_{\eps}(x)=\phi_0(x)+\eps \varphi(x)$, where the function $\varphi$ is given by Proposition~\ref{p:non-analytic}. Then at any point $(x,\alpha_0)$ for all values of the parameter $\eps$ except at most one, the corresponding function $\Psi(\phi_{\eps})$ is non-analytic in~$\alpha$.
\end{coro}
\begin{proof}
    Indeed, if the $\Psi$-images of two of the functions $\phi_{\eps_1}$, $\phi_{\eps_2}$ were analytic in $\alpha$ at some point $(x,\alpha_0)$, then the image of their difference would be also analytic in $\alpha$ at the same point, contradicting the conclusion of Proposition~\ref{p:non-analytic}.
\end{proof}
\begin{proof}[Proof of Proposition~\ref{p:non-analytic}]
    Note that the operator $\Phi$ commutes with rotations for every fixed $\alpha$, and thus it should be diagonal in the Fourier basis. 
     For simplicity of the formulae (so that we can use complex exponents), let us temporarily allow the function $\varphi$ to be complex-valued.
     %; passing to its real or imaginary part would conclude the proof for the real-valued functions, too.
    Writing explicitly the images of a function 
    $e^{2\pi i m x}$, one has  %\varphi_m:=
    \[
        \Psi(e^{2\pi i m x}) = \sum_{n=1}^{\infty} \frac{e^{-2\pi i m \alpha \cdot n}}{2^n} \, e^{2\pi i m x} = \frac{\frac{1}{2}e^{- 
        2\pi i m \alpha}}{1- \frac{1}{2}e^{- 
        2\pi i m \alpha}} \, e^{2\pi i m x}.
    \]
     Due to this explicit form as a function of $\alpha$, one can note that for every $m$ and every $x_0$ the function $\Psi(e^{2\pi i mx}) (x_0,\cdot)$ has poles at the points of the set
    \[
        A_m = \{ \alpha \mid e^{-2\pi i m \alpha} = 2\} 
        = \frac{\ln 2}{2\pi m} i + \frac{1}{m} \Z.
    \]
    Returning to the real-valued functions, we have 
    \[
        \Psi(\cos 2\pi m x) = \Psi(\frac{1}{2}e^{2\pi i m x}) + \Psi(\frac{1}{2}e^{-2\pi i m x}),
    \]
    so the poles $\Psi(\cos 2\pi m x)$ are exactly the points of the set 
    \[
        A'_m= A_m \cup \overline{A_m}.
    \]
    
%    \missingfigure{Poles $A_m$ and the distances}
    
    Now, the radius of convergence of a holomorphic function is given by the distance to its closest singularity, and for any real $\alpha_0$ its distance to the set of poles $A'_m$ one has (see Fig.~\ref{fig:Poles})
    \[
       \frac{1}{10m}<\frac{\ln 2}{2\pi m}\le \dist(\alpha_0,A'_m) \le \frac{\sqrt{1+(\frac{\ln 2}{\pi})^2}}{2m} <\frac{0.6}{m}. 
    \]
    \begin{figure}
        \centering
        \includegraphics[width=0.8\linewidth]{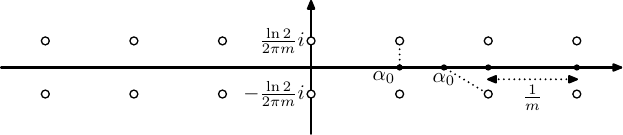}
        \caption{The set of poles $A'_m$ (empty circles) and the points $\alpha_0\in \R$, maximizing and minimizing the distance to it.}
        \label{fig:Poles}
    \end{figure}
    
    % In particular, for any $m$
    % \[
    %     \forall \alpha_0,x \quad  \forall N \quad \exists n>N: \quad |\Ta_n(e^{2\pi i m x},x,\alpha_0)| > \left(\frac{m}{\sqrt{2}}\right)^n.
    % \]
    % Moreover, due to the continuity of the functions and compactness of the circle, one has
    % \[
    %     \forall N \quad \exists N': \forall \alpha_0,x \quad \exists n, \, N<n<N': \quad |\Ta_n(e^{2\pi i m x},x,\alpha_0)| > \left(\frac{m}{2}\right)^n.
    % \]

     We are going to recurrently construct a sequence $(b_j)_{j\in \N}$ of Fourier coefficients, associated to harmonics with numbers $m_j:=10^j$, such that
  for the function $\varphi$, given by the Fourier series 
    \begin{equation}\label{eq:vphi-def}
        \varphi(x):=\sum_{j=1}^{\infty} b_j \cos 2\pi m_{j}x,
    \end{equation}
  the Taylor coefficients in $\alpha$ of its image $\Psi(\phi)$ grow superexponentially at every~$(\alpha_0,x)$. To guarantee such growth, we will at the same time construct a sequence $(N_j)_{j=0,1,\dots}$ of degrees, and ensure that a ``large'' Taylor coefficient appears somewhere between degrees~$N_{j-1}$ and~$N_j$:
  \begin{equation}\label{eq:growth-between-N}
    \forall \alpha_0,x\in \R \quad \forall j \quad  \exists k, \, N_{j-1}<k<N_j: \quad |\Ta_k(\varphi,x,\alpha_0)|> (m_j)^k.      
  \end{equation}
  Clearly, once~\eqref{eq:growth-between-N} is established, it implies the desired~\eqref{eq:c-growth} immediately, thus concluding the proof.
  %More specifically, we will also construct a sequence of degrees $(N_j)_{j=0,1,\dots}$, and ensure with each $b_j$ that a ``large'' Taylor coefficient appears somewhere between degrees~$N_{j-1}$ and~$N_j$. Formally speaking, the choice is made to satisfy the properties~\eqref{eq:Taylor-lower} and~\eqref{eq:Taylor-upper} below, that will guarantee the conclusion of Proposition~\ref{p:non-analytic} for the function
%    where $m_j:=2^j$.
%\todo[inline]{Polish the text below, check constants} 

    To perform such a construction, we start by setting $N_0=0$, $b_1=1$. Then,  for every $j=1,2,\dots$ we recurrently:
    \begin{itemize}
        \item Choose $N_j>N_{j-1}$ sufficiently large so that 
        \begin{multline}\label{eq:Taylor-lower}
            \forall \alpha_0,x \quad \exists k, \, N_{j-1}<k<N_j: \quad \left|\Ta_k\left(\sum_{r\le j} b_r \cos{2\pi m_r x},x,\alpha_0\right)\right| \ge \\  
            b_j \cdot \left|\Ta_k\left(\cos {2\pi m_j x},x,\alpha_0\right)\right| - \sum_{r\le j-1} b_r \cdot \underbrace{\left|\Ta_k\left(\cos {2\pi m_r x},x,\alpha_0\right)\right|}_{\displaystyle \preccurlyeq
            (10 m_r)^k}> 2\cdot \left(m_j\right)^k;           
        \end{multline}

        % \begin{equation}\label{eq:Taylor-lower}
        %     \forall \alpha_0,x \quad \exists n, \, N_{j-1}<n<N_j: \quad \left|\Ta_n\left(\sum_{r\le j} b_r e^{2\pi i m_r x},x,\alpha_0\right)\right| > \left(\sqrt{2}\pi m_j\right)^n;     
        % \end{equation}
        \item Then choose $b_{j+1}\in (0,2^{-m_{j}})$ sufficiently small so that 
        \begin{equation}\label{eq:Taylor-upper}
            \forall \alpha_0,x \quad \forall n<N_{j} \quad | \Ta_n(b_{j+1} \cos {2\pi m_{j+1} x},x,\alpha_0)|< \frac{1}{2^j}.        
        \end{equation}
    \end{itemize}
    The fact that we can satisfy~\eqref{eq:Taylor-lower} follows from joining the upper estimate by $\frac{0.6}{m_j}$ of the radius of convergence of $\Psi(\cos {2\pi m_j x})$ with a lower estimate by $\frac{\ln 2}{2\pi m_{r}} > \frac{1}{m_j}$ for the radii of convergence for all $\Psi(\cos {2\pi m_r x})$ with $r=1,\dots,j-1$. The fact that we can satisfy~\eqref{eq:Taylor-upper} is immediate due to the linearity on~$b_{j+1}$: 
    \[
        \Ta_n(b_{j+1} \cos {2\pi m_{j+1} x},x,\alpha_0) = b_{j+1} \Ta_n(\cos {2\pi m_{j+1} x},x,\alpha_0),
    \]
    so~\eqref{eq:Taylor-upper} is equivalent to 
    \begin{equation}\label{eq:b-j}
        b_{j+1} \cdot \max_{n\le N_j} \sup_{x,\alpha_0\in \Sc} |\Ta_n(\cos {2\pi m_{j+1} x},x,\alpha_0)| <\frac{1}{2^j},
    \end{equation}
    and in the left hand side of~\eqref{eq:b-j} we have a maximum of supremums a finite number of continuous functions, considered over a compact set $\Sc\times \Sc$, that is thus finite. 

    Now, take $\varphi(x)$ to be defined by~\eqref{eq:vphi-def}.
    This function is analytic, as its Fourier coefficients decay exponentially.
    At the same time, joining~\eqref{eq:Taylor-lower} and~\eqref{eq:Taylor-upper}, we get 
    \begin{multline*}
        \forall x,\alpha_0 \quad \forall j \quad \exists k, \, N_{j-1}<k<N_j : \quad |\Ta_k(\varphi,x,\alpha_0)|> 
        \\
        2 \left(m_j\right)^k - \sum_{r>j} \frac{1}{2^r} > 2 \left(m_j\right)^k - 1> m_j^k ,
    \end{multline*}
    thus implying the desired~\eqref{eq:growth-between-N} and thus concluding the proof of~\eqref{eq:c-growth}.
\end{proof}

 %{Finally, similar arguments can be repeated with functions~$\cos (2\pi m x)$ used instead of $e^{2 \pi i m x}$, providing a real-valued function~$\varphi$.} 
%or $\sin (2\pi m x)$ 

    % that the Taylor where the index $N_j$ guarantees that the Taylor coefficient, associated to $b_j e^{2\pi i m_j x}$, reaches a sufficiently high value between $N_{j-1}$ and $N_j$,
    % % \begin{equation}\label{eq:Taylor-lower}
    % %     \forall m \quad \forall \alpha_0,x \quad \exists n, \, N_{m-1}<n<N_m: \quad | \Ta_n(b_m e^{2\pi i m x},x,\alpha_0)| > \left(\frac{m}{2}\right)^n,        
    % % \end{equation}
    % while $b_m<\frac{1}{2^m}$ and its effect on smaller indices is sufficiently small,
    % \begin{equation}\label{eq:Taylor-upper}
    %     \forall m \quad \forall \alpha_0,x \quad \forall n<N_{m-1} \quad | \Ta_n(b_m e^{2\pi i m x},x,\alpha_0)|< \frac{1}{2^m}.        
    % \end{equation}
    % Now, define $\varphi(x)=\sum_m b_m e^{2\pi i m x}$. The Fourier coefficients $b_m$ decrease exponentially, and hence the function $\varphi(x)$ is analytic. At the same time, for each $m$

\begin{remark}
     It seems highly plausible that actually the function $\Psi(\varphi)$ is non-analytic in $\alpha$ for all~$\varphi$, except for the trigonometric polynomials (due to the accumulation of singularities of the $\Psi$-images of Fourier harmonics $e^{2\pi i m x}$ to the $\alpha$-real axis). Though, proving such a statement rigorously might be technically complicated.
\end{remark}

\section{Concluding remarks}\label{s:conclusion}

Here we provide a few concluding remarks and open questions.

\subsection{Dynamical arguments and related questions}

    The limits at rational frequencies in Theorem~\ref{t.Leb-cont} in the subcritical regime $\lambda<1$ can be considered from the dynamical viewpoint on the corresponding cocycle. As the simplest example, consider the limit $\alpha\to 0$. Then, if some energy $E$ does not belong to the intersection spectrum, there is an interval $I\subset S^1$, such that for phases $\theta\in I$ the corresponding matrices $\Pi_{E,\theta}$ are hyperbolic.  Then, for the transfer matrix along the period $\alpha=\frac{p}{q}, \, p\ll q,$ among the multiplied matrices $\Pi_{E,j\alpha+\theta}$ the proportion of the hyperbolic ones is at least~$|I|$. For such a product to be elliptic, the energy $E$ should be very finely tuned, so that the elliptic part of the product would send the expanded direction of the hyperbolic part of the product close to the contracted one. 
    
     This is only a hand-waving description; for instance, handling frequencies of the form $\frac{p}{np+1}$, for which the orbits make $p$ turns around the circle, would be already a nontrivial task. However, we hope that it explains why the intervals of the spectrum $\Sigma_{\alpha,\lambda}$ away from the intersection spectrum at
    %, say %$\alpha_0=0$,
    zero frequency are very small (and hence carry a Lebesgue measure that tends to zero).

% In the subcritical regime, Chambers' formula implies that traces, associated to the long periodic products, almost do not depend on the phase (due to exponentially small term~$\lambda^q$ in~\eqref{eq:Chambers}), and this helps with the above hand-waving explanation (hyperbolicity more depends on the phas
% Contrary to the above  explanation, 
In the subcritical regime, due to exponentially small term~$\lambda^q$ in~\eqref{eq:Chambers}, the hyperbolicity or ellipticity of the long periodic product almost does not depend on the initial phase (and that contributes to the above explanation).
The situation in the supercritical regime is drastically different, and we have no even hand-waving dynamical argument for the support of the corresponding $\mu_{\frac{p}{q},\lambda}^-$. 

\begin{question}
    In the supercritical regime $\lambda>1$, for $\alpha=0$, the support of the measure 
    $\mmu_{0,\lambda}$ is the interval $[-2(\lambda-1),2(\lambda-1)]$. Is there a good dynamical explanation for that? Is it related to the fact that for these energies the ``hyperbolic domain'' on the phase circle consists of two disjoint arcs, separated by ``elliptic'' arcs?
\end{question}

Also, note that the arguments of the present paper rely on the explicit form of the Almost Mathieu potential, $V(\theta)=2\lambda\cos\theta$: Remark~\ref{r:Si-Q} shows that one of the starting arguments fails even for the potentials that are linear combination of $\cos 2\pi \theta$ and $\cos 4\pi \theta$. In the meantime, dynamical (cocycles) point of view theoretically might be applicable for other potentials as well. We thank S. Jitomiskaya for suggesting the following question.
\begin{question}
    Does Theorem~\ref{t.Leb-cont} above generalise to other subcritical potentials? For instance, to $V(\theta)$ that are (even) trigonometric polynomials? Any analytic potential?
\end{question}

\subsection{The explicit form of the moments of $\mmu_{\alpha,\lambda}$}

Continuing the sequence of moments of the measures $\mmu_{\alpha,\lambda}$ after the second and the fourth ones, discussed in Example~\ref{ex:2-4}, one gets: 
\begin{example}\label{ex:6-8}
    The sixth and eights moments in the subcritical regime $\lambda<1$ are given by 
    \[
        \frac{2k+1}{2^{2k+2}} c_{2k}(\alpha,\lambda) = T_{2k}(\alpha,\lambda)-\lambda^{2k+1} T_{2k}(\alpha,\frac{1}{\lambda}), \quad k=3,4, %,5
    \]
    where, using the notation $t:=\cos 2\pi \alpha$, one has
    \begin{multline}\label{eq:T-6}
        T_6(\alpha,\lambda)
                =1+ \frac{7}{2}\left(t^5+t^4+t^3+t^2+t+1 \right) \lambda^2 
        \\ +\frac{7}{4} (1+t)^2 \cdot \left( 4t^4 -4 t^3 + 3t^2 + 2 \right) \cdot \lambda^4 + \frac{7}{8} (1+t)^3 \lambda^6,
    \end{multline}
    and 
    \begin{multline}\label{eq:T-8}
        T_8(\alpha,\lambda)
                =1+ \frac{9}{2}\left(t^7+ t^6+ t^5+t^4+t^3+t^2+t+1\right) \lambda^2 
                \\ + \frac{9}{4} (1+t)^2 \left( 16 t^8- 16 t^7- 4 t^6 + 8 t^5+9 t^4-8 t^3+6 t^2+3 \right) \lambda^4 
        \\ + 9 (1+t)^3 \cdot \left( 2 t^6-4 t^5+4 t^4-2 t^3+ \frac{3}{4}t^2+ \frac{5}{12} \right) \cdot \lambda^6 + \frac{9}{16} (1+t)^4 \lambda^8,
    \end{multline}    
    % and 
    % \begin{multline}\label{eq:T-10}
    %     T_{10}(\alpha,\lambda)
    %             =1+ \frac{11}{2}\left(t^9+ t^8+ t^7+ t^6+ t^5+t^4+t^3+t^2+t+1\right) \lambda^2 
    %             \\ + \frac{11}{4} (1+t)^2 \Bigl( 64t^{12}-64t^{11}-80t^{10}+96t^9+68t^8
    %             \\ -92t^7-5t^6+32t^5+10t^4-12t^3+9t^2+4 \Bigr) \lambda^4 
    %             \\ + \frac{11}{8} (1+t)^3 \Bigl( 256t^{12}-512t^{11}+128t^{10}+320t^9-112t^8
    %             \\ -128t^7+68t^6-32t^5+75t^4-48t^3+20t^2+7\Bigr) \lambda^6 
    %     \\ + \frac{11}{2^4} (1+t)^4 \cdot \left( 64t^8-192t^7+272t^6-224t^5+120t^4-40t^3+10t^2+5 \right) \cdot \lambda^8 
    %     \\ + \frac{11}{2^5} (1+t)^5 \lambda^{10}.
    % \end{multline}    
\end{example}
    Let us mention some initial %There are some first 
    observations that can be made looking at % from
    these examples. First, one notices that $\lambda^2$ in $T_{2k}(\alpha,\lambda)$ comes with the factor $\frac{2k+1}{2} \cdot \frac{1-t^{2k}}{1-t}$ and that $\lambda^{2k}$ comes with the factor $\frac{2k+1}{2^k} (1+t)^k$. A way to prove the latter is by using that
    \[
           T_{2k}(\alpha,\lambda)=\frac{1}{2^{2k+1}} \tr H_{\alpha,\lambda,\frac{\alpha}{2}}^{2k+1} R_{\frac{\alpha}{2}}, 
    \]
    and noticing that for $\lambda^{2k}$ terms it suffices to consider the submatrix with indices $i=-1,0$. Hence, we are looking for the coefficient at $\lambda^{2k}$ of the polynomial % ; one thus gets
    \begin{multline}        
         \frac{1}{2^{2k+1}} \tr 
            \begin{pmatrix}
           2\lambda\cos \pi \alpha & 1 
            \\
           1 & 2\lambda\cos \pi \alpha
       \end{pmatrix}^{2k+1} 
        \begin{pmatrix}
           0 & 1 
            \\
           1 & 0
       \end{pmatrix} 
              \\
       = \frac{1}{2^{2k+1}}  \left((2\lambda \cos \pi \alpha +1)^{2k+1}- (2\lambda \cos \pi \alpha -1)^{2k+1}\right),
       \end{multline}
       which is equal to 
    $$
        \frac{1}{2^{2k+1}} \cdot 2 \cdot (2k+1) (2\cos \pi \alpha)^{2k} = \frac{2k+1}{2^k} (1+t)^k.
   $$ 
    For the  coefficient at $\lambda^2$, it can be found by applying the generalized Aubry--Andr\'e duality, as it could be seen as the $\lambda^{2k-1}$ term in
    \[
        \lambda^{2k+1} T_{2k}\left(\alpha,\frac{1}{\lambda}\right) = \frac{1}{2^{2k+1}} \tr H_{\alpha,\lambda,0}^{2k+1} R_{0},
    \]
    and to account for such terms, it suffices to consider the submatrix with indices $i=-1,0,1$, which gives a polynomial
    \[
              \frac{1}{2^{2k+1}} \tr 
            \left[\begin{pmatrix}
           2\lambda t & 1 & 0
            \\
           1 & 2\lambda & 1
            \\
           0 & 1 & 2\lambda t
       \end{pmatrix}^{2k+1} 
        \begin{pmatrix}
           0 & 0 & 1 
            \\
           0 & 1 & 0 
            \\
           1 & 0 & 0
       \end{pmatrix}\right] 
            \] 
    with a coefficient at $\lambda^{2k-1}$ equal to 
    $$
    \frac{2k+1}{2} (1+t+\dots+t^{2k}).
    $$
    In a similar way one could theoretically obtain expressions for coefficients for instance, for $\lambda^4$ or $\lambda^{2k-2}$, but the computations quickly become exceedingly technical. 
    
    A next observation is that the $\lambda^{2j}$ term is divisible by $(1+t)^j$; it can be again seen from interpreting $T_{2k}(\alpha,\lambda)$ as $\frac{1}{2^{2k+1}} \tr H_{\alpha,\lambda,\frac{\alpha}{2}}^{2k+1} R_{\frac{\alpha}{2}}$. Indeed, the diagonal elements of the matrix $H_{\alpha,\lambda,\frac{\alpha}{2}}$ are of the form $2\lambda \cos (2m+1) \pi \alpha$, and $\cos (2m+1) \pi\alpha$ is given by an odd polynomial of $\cos \pi \alpha$, hence the latter factors out. This leads to $\lambda^{2j}$ term having a factor $\cos^{2j} \pi\alpha = \left(\frac{1+t}{2}\right)^j$.

\begin{question}
    Can any interesting conclusions be made about the polynomials that are left after factoring $(1+t)^{2j}$ from $\lambda^{2j}$ term in $T_{2k}$? For instance, they all seem to be missing a linear term in $t$; is there anything else that can be said? Is there any direct way to produce the polynomials $T_{2k}$, either by an explicit formula or via some recurrence relation?
\end{question}

\section*{Acknowledgments}
The authors are grateful to Lana Jitomirskaya for providing multiple useful references.

\end{document}